\documentclass[reqno,12pt]{article}
\usepackage{amsmath,amsthm,amsfonts,amssymb,amscd}
\usepackage[dvips]{graphicx,color,epsfig}
\usepackage[T1]{fontenc}
\usepackage[latin1]{inputenc}
\newtheorem{thm}{Th\'eor\`eme}[section]
\newtheorem{cor}[thm]{Corollaire}

\newtheorem{lem}[thm]{Lemme}

\newtheorem{pro}[thm]{Proposition}

\newtheorem{defn}[thm]{D\'efinition}

\newtheorem{rem}[thm]{Remarque}
\newtheorem{rems&defn}[thm]{Remarque et d\'efinition}

\newcommand{\R}{\mathbb{R}}
\newcommand{\EE}{\mathbb{E}}
 \newcommand{\LL}{\mathbb{L}}
\newcommand{\F}{\mathcal{F}}

\newcommand{\A}{{\cal A}}
\newcommand{\T}{\mathcal{T}}
\newcommand{\TT}{\mathcal{T}_t}

\newcommand{\G}{\mathcal{G}}
\newcommand{\ESS}[1]{\underset{#1}{\mathrm{essup}}}
\newcommand{\INF}[1]{\underset{#1}{\mathrm{ess~inf}}}

\newcommand{\demi}{\frac{1}{2}}
\newcommand{\eps}{\varepsilon}
\newcommand{\1}{{\mathbf 1}}

\newcommand{\To}{\longrightarrow}

\newcommand{\biindice}[3]
{ #1\begin{array}[t]{c}
{\scriptstyle #2}\\
{\scriptstyle #3}\end{array}
}

\def \endproof {\quad \hfill  \rule{2mm}{2mm}\par\medskip}
\begin{document}

\title{Innovation, commutation et contrôle impulsionnel en horizon infini}
\author{Rim Amami \\
Institut de Mathématiques de Toulouse \\
Université Paul Sabatier,  Toulouse\\
rim.amami@math.univ-toulouse.fr }
\maketitle

\begin{abstract}

 Nous étudions  un problème de contrôle impulsionnel en horizon infini. 
Pour le résoudre, nous étendons au cas de l'horizon infini des résultats
concernant  les équations différentielles stochastiques rétrogrades et 
réfléchies à double barrière.  Les propriétés de l'enveloppe de Snell permettent  de réduire notre problème
 à montrer l'existence d'un couple de processus continus. 
\end{abstract}
   \textbf{ Mots clés:} Contrôle impulsionnel, horizon infini,   équations différentielles stochastiques      rétrogrades et réfléchies,  double barrière.

\section{Introduction}

La théorie du contrôle impulsionnel est associée au problème du choix optimal de temps successifs 
pour effectuer une action (dans notre cadre l'action est le changement de technologie)  afin d'optimiser un critère économique. C'est une situation où l'on fait face à des systèmes dynamiques évoluant dans des conditions d'incertitude et où il s'agit de prendre des décisions à chaque date afin de maximiser une espérance de gain. Ce type de problème se résout entre autres en utilisant l'enveloppe de Snell et le principe de programmation dynamique.
%

Le problème de contrôle impulsionnel est un sujet qui apparaît souvent dans la littérature par exemple  l'économie, la statistique et  les mathématiques financières.  Il a été  initi\'e par  Bensoussan et  Lions
\cite{bensoussan} et ensuite formalisé par d'autres auteurs. 

  Parmi d'autres,  Hamadène et al. \cite{hamadene}   ont montré l'existence d'une solution     pour les équations différentielles stochastiques rétrogrades   réfléchies  à double barrière en horizon fini en utilisant les propriétés de l'enveloppe de Snell.

 Hdhiri et  Karouf \cite{monia} ont prouvé  l'existence d'une stratégie optimale maximisant le gain de la firme dans un cadre non-markovien en utilisant une fonction d'utilité de type exponentielle et des propriétés de l'enveloppe de Snell.  
  
  D'autres auteurs ont utilisé les équations stochastiques différentielles rétrogrades  pour leurs modèles. Par exemple, 
   Jeanblanc et  Hamadène    ont considéré  dans  \cite{ jeanblanc} une centrale électrique  qui a deux modes: arrêt et démarrage. Il s'agit d'un problème de contrôle impulsionnel avec changement du contrôle sans saut de la variable d'état. 
  Ce problème a été  résolu en utilisant principalement l'enveloppe de Snell et les équations différentielles stochastiques   réfléchies.
  
       Cvitanic et    Karatzas \cite{cvitanic } ont établi l'existence et l'unicité des solutions des équations différentielles stochastiques  rétrogrades à deux barrières réfléchies.  Ces auteurs ont généralisé les travaux de  El Karoui et al. \cite{Karoui1}. Ils ont aussi montré que la solution des EDS coincide avec la valeur du jeu de Dynkin, un jeu stochastique d'arrêt optimal. 
  
   Djehiche et al.  \cite{djehiche}   ont 
 utilisé des outils purement probabilistes comme les équations différentielles stochastiques rétrogrades et l'enveloppe de
Snell     pour
résoudre le problème optimal de changement de technologie en horizon fini.  
 
 Mnif et al. \cite{mnif2} ont utilisé le principe de programmation dynamique  des équations  Hamilton-Jacobi-Bellman. Ces auteurs ont fourni un    schéma  numérique pour les 
inéquations quasi variationnelles  associées aux problèmes de contrôle   impulsionnel.

Vath et   Mnif  ont prouvé dans \cite{mnif1} que le problème de contrôle   impulsionnel est réduit à une suite itérative des problèmes d'arrêt optimaux. En effet, en utilisant des inégalités variationnelles, la fonction de valeur est obtenue comme la limite d'une suite itérative de ces problèmes d'arrêt optimaux. Ensuite, les auteurs ont résolu  le problème d'arrêt numériquement à l'aide de la méthode de Monte Carlo et le calcul de Malliavin. 

  Bouchard et  Chassagneux \cite{bouchard} ont étudié l'approximation en temps discret de la solution d'une équation différentielle stochastique Forward Backward réfléchie dans le cas où la réflexion n'opère que dans un ensemble fini de fois. Les auteurs ont prouvé des propriétés de régularité de la solution de l'EDS sous des hypothèses de Lipschitz. D'autres auteurs ont étudié l'approximation en temps discret  des EDS Forward Backward réfléchies, citons parmi d'autres \cite{bouchard1, bouchard2, chassagneux}.  
      
 Les  équations différentielles stochastiques  rétrogrades (EDSR) ont été   initiées   par   Pardoux et Peng \cite[année 1990]{pardoux}, \cite[année 1992]{pardoux1}. Ils  ont été 
  les premiers à résoudre   le problème  de  l'existence et l'unicité de la solution des EDSR sous des hypothèses de    Lipschitz sur la fonction drift. 
  Plusieurs auteurs ont été  attirés par ce domaine notamment dans le cadre des applications numériques,  notamment pour  améliorer les conditions d'existence et d'unicité d'une solution de l'EDSR. Tous ces travaux ont été  basés principalement sur le théorème de comparaison des solutions des EDSR. Parmi d'autres, 
  citons \cite{Karoui2, peng, yong}.  
  
   El Karoui et al. \cite{Karoui1} et  Hamadène et al. \cite{hamadene} ont étendu ces résultats aux cas des EDS  rétrogrades  et réfléchies.   Néanmoins,  leurs résultats ne s'appliquent pas  directement à la  situation qui nous intéresse qui exige un horizon infini.   Ceci dit, leurs 
    papiers ont  fourni beaucoup d'inspiration et de motivation à  notre travail. Ces auteurs  fournissent  une solution au problème   des EDS rétrogrades et réfléchies  que nous étendons au cas de l'horizon infini   en ajoutant un coefficient d'actualisation et en imposant des conditions d'admissibilité sur  les stratégies, c'est à dire des suites de temps d'impulsion. \\
   
        Notre  principale contribution   est de prouver l'existence  d'une stratégie optimale 
  qui maximise le gain moyen d'une  firme mais, à la différence des précédents auteurs, 
  dans un contexe en horizon infini.

 Notre travail est organisé comme suit:  la  section \ref{model}  est
consacrée à définir le modèle de contrôle impulsionnel à résoudre.  Dans  la section \ref{snell},  nous rappelons quelques notions fondamentales de l'enveloppe de Snell et     nous   présentons   un couple de processus continus dont l'existence, prouvée en section \ref{existe}, permettra
d'exhiber une stratégie optimale.  
   Dans la section \ref{EDSFB}, nous   généralisons les outils des équations rétrogrades  et  réfléchies à double barrière à l'horizon infini sous des hypothèses convenables.
  
  \section{Présentation du  modèle}
  \label{model}
  Soit $(\Omega,\F, (\F_t)_{t \geq
0}, \mathbb{P})$ un espace de probabilit\'e
 muni d'une filtration $(\F_t)_{t \geq
0}$ compl\`ete continue \`a droite  et soit un mouvement brownien $W
= (W_t)_{t \geq 0 }$.  
Nous noterons par $(\G_t)_{ t  > 0}$ la filtration
 d\'efinie par $\G_t = \displaystyle\vee_{s < t}\F_s,\; \forall\; s<t$, (notée usuellement $\F_{t^-}$).\\
 \\
 Supposons qu'une  entreprise  décide à des  instants aléatoires de changer de technologie afin de maximiser son gain.   Nous supposons que $\{0,1\}$ (noté $U$)
 est l'ensemble des technologies permises telles que $0$ est l'ancienne technologie et $1$ est la nouvelle technologie. 
  L'évolution de la firme dépendant de plusieurs facteurs externes (prix sur le marché, crise mondiale, temps,...), 
  le changement de technologie provoque un coût défini par $c_{0,1}$ si on passe de la technologie $0$ à la technologie $1$ et $c_{1,0}$ 
  si on passe de la technologie $1$ à la technologie $0.$ Supposons que $c_{0,1} > c_{1,0}.$  
  
 Nous définissons  une stratégie de contrôle impulsionnel  comme une suite:
$$\alpha := (\tau_n)_{n\geq -1},$$
 où  $(\tau_n)_{n \geq -1}$ est une suite croissante de
$\G$-temps d'arrêt  avec   $\tau_{-1} = 0$. On note 
$$\tau:= \lim_{n \To +\infty} \tau_n.$$ 
 La suite $(\tau_n)$  
 modélise  la suite des instants d'impulsion du système de la façon suivante:
   pour tout $n \geq 0,$    $\tau_{2 n }$ est l'instant où la firme passe de la technologie $0$ à la technologie $1$ et  $\tau_{2 n +1}$ est  l'instant   où la firme passe de  $1$ à   $0$.  
   
     Nous introduisons le processus
c\`adl\`ag $\xi$ à valeurs dans  $\{0,1\}$ défini par    $\xi_0= 0$  et

\begin{equation}
\label{xit}
  \xi_t := \xi_0 \1_{[0, \tau_{0} [} (t)+  \sum_{n\geq 0} \1_{[ \tau_{2n}, \tau_{2n+1} [}  (t)  + \xi_{\tau^-} \1_{[\tau,+\infty[}(t).
\end{equation}
La valeur de la firme     est donnée par  
 $S_t = \exp X_t^\xi,~ t\geq 0,$
où $X^\xi$ est  le processus continu  à droite défini par:
\begin{equation}\label{Y}
     d X_t^\xi  = b(\xi_t,X^{\xi} _t)\,dt + \sigma(\xi_t,X^{\xi} _t)\,dW_t,
\end{equation}
où   $b:{U}\times \mathbb{R}\rightarrow \mathbb{R}$ et $\sigma:
{U}\times \mathbb{R}\rightarrow \mathbb{R}^{+}$  sont deux fonctions   
mesurables satisfaisant  sur $\R$ la condition  de Lipschitz et  la condition   de croissance sous-linéaire:\\
  - Il existe une constante   $  K \geq 0$ telle que pour tout $i  \in {U}$ et tout  $  x,y \in \mathbb{R}$,
$$
\big|b(i,x) -  b(i,y) \big|  + \big|\sigma(i,x) -  \sigma(i,y) \big|\leq K \big|x-y\big|.
$$
 -  Il existe une constante   $  K \geq 0$ telle que pour tout $i  \in {U}$ et tout  $  x  \in \mathbb{R}$,
$$
\big|b(i,x)\big|^2  + \big|\sigma(i,x)   \big|^2 \leq K^2 (1+ |x |^2). 
$$ 

En appelant     $f>0$ le bénéfice net de la firme et    $c>0$   le coût de changement de technologie,      toute strat\'egie $\alpha$ occasionne un gain:
 $$
  \int_0^{+\infty}   e^{ -\beta s} f(\xi_s,X_s^{\xi_s} ) ds -
   \sum_{n\geq 0}  \left\{  e^{ -\beta \tau_{2n}} c_{0,1}     +
    e^{ -\beta \tau_{2n+1 }} c_{1,0}     \right\},
$$
où $\beta>0$ est un coefficient d'actualisation.
\begin{defn}\label{strategadm}
  La stratégie $\alpha = (\tau_n)_{n\geq -1} $   est dite admissible si et seulement si: 
 $$ 
  \mathbb{E} \int_0^ {+\infty} 
   e^{-\beta s}  f(\xi_s,X_s^{\xi_s} )  ds  < \infty \qquad  \mathbb{E} \sum_{n\geq 0} \{    e^{ -\beta \tau_{2n}} c_{0,1}     +
    e^{ -\beta \tau_{2n+1 }} c_{1,0}  \}  < \infty.
$$
On note $\A$ l'ensemble des strat\'egies
admissibles. 
\end{defn}

Le problème de contrôle impulsionnel  posé consiste à prouver l'existence
d'une stratégie admissible  $\widehat{\alpha}$ qui maximise la
fonction gain  moyen $ K(\alpha,i,x)$ d\'efinie par
\begin{equation}
\label{Kalpha}
  K(\alpha,i,x ) = \ESS{\alpha \in \A} ~ \mathbb{E}_{i,x}\left[   \int_0^{+\infty}   e^{ -\beta s} f(\xi_s,X_s^{\xi_s}  ) ds -
   \sum_{n\geq 0}  \left\{  e^{ -\beta \tau_{2n}} c_{0,1}     +
    e^{ -\beta \tau_{2n+1 }} c_{1,0}     \right\}     \right].
   \end{equation}
   \\
 Notons que sur l'événement $\{\tau<\infty\}$ le gain après $\tau$
est $$ \mathbb{E}_{i,x}[\int_\tau^{+\infty}   e^{ -\beta s} f(\xi_s,X_s^{\xi_s}  ) ds  \big|  \F_{\tau}]:$$
c'est le résultat de la décision étant de ne plus changer de stratégie.  \\

  Enfin, introduisons  les ensembles utiles suivants:
  \begin{itemize}
    \item [$\bullet$] $\T= \big\{ \theta: \G\mbox{-temps d'arrêt} \big \}.$
    \item [$\bullet$]  $\TT= \big\{ \theta\in \T~:  \theta \geq t  \big \}.$
       \item [$\bullet$] $\mathcal{P}^2= \{\mbox{  processus} ~ \F_.\mbox{-progressivement  mesurables } \}$
    \item [$\bullet$] $\mathcal{C}^2= \big\{  (X_t)_{t \geq 0} \in \mathcal{P}^2 :  ~ \mbox{  tel que } ~   \mathbb{E}[ \sup_{t \geq 0} |X_t|^2] < \infty  \big \}.$
   \item [$\bullet$] $\mathbb{H}^2= \big\{  (X_t) \in \mathcal{P}^2:   \mbox{  tel que}  ~   \mathbb{E}[ \int_0^{\infty }  |X_t|^2 dt ] < \infty  \big \}.$
   \item [$\bullet$] $\mathcal{D}_n= \big\{  \mu \in \mathcal{P}^2: [0,T]\times\Omega \To [0,n] \big \}.$    
\end{itemize}

  \section{Enveloppe de Snell}
  \label{snell}
  
 Commençons par rappeler quelques notions fondamentales de contrôle optimal de   El Karoui  \cite{Karoui} utilisées pour résoudre les problèmes de contrôle impulsionnel étudiés dans ce chapitre.
 
 \begin{defn} 
 Un processus  $U$ est  de classe   (D)   si l'ensemble des  variables aléatoires  $\{U_{\theta}, \theta \in \mathcal{T}\}$ est uniformément intégrable.
 \end{defn}
 
 \begin{thm}  Soit   $U$ un processus $\F$-adapté, càdlàg de classe (D). Notons $Z$ son  enveloppe de Snell. C'est  la plus petite sur-martingale de classe (D) qui majore  $U$:
 \begin{equation}
 \label{envsnell}
 Z_t= \ESS{\theta \in  \TT}  ~  \mathbb{E}[U_{\theta} \big| ~ \F_t].
\end{equation}
\end{thm}

  \begin{rem}\label{decomp}  Le processus $Z$ admet la décomposition unique suivante:
  $$Z = M - A,$$
  où:
    \begin{itemize}
    \item [(i)] $M$ est une martingale.
    \item [(ii)] $A$ est un processus croissant, continu à droite, intégrable et $A_0=0$.
         \end{itemize}
  \end{rem}
  
\begin{defn}    \label{defopt}
 Soit U un processus $\F$-adapté et $Z$ son    enveloppe de Snell.  Une condition nécessaire et suffisante pour qu'un temps d'arrêt $\widehat{\tau}$
             soit optimal après $\gamma $  est que:
     \begin{itemize}
    \item [$\bullet$] $\widehat{\tau} \geq \gamma, \quad  \gamma \in  \mathcal{T}$. 
    \item [$\bullet$]  $Z_{\gamma}= \mathbb{E}[Z_{\widehat{\tau} } \big| ~ \F_{\gamma}]=    \mathbb{E}[U_{\widehat{\tau} } \big| ~ \F_{\gamma}].$
             \end{itemize}
             En particulier $Z_{0}=   \sup_{\theta \in \TT } \mathbb{E}[U_{\theta}  ] = \mathbb{E}[U_{\widehat{\tau} }].$ Par exemple, 
              $\widehat{\tau} 
             = \inf \{ t \geq \gamma: Z_t = U_t\} $ 
             est le premier temps optimal après $\gamma.$
           \end{defn}
           
Le but principal de ce chapitre est de prouver l'existence d'une stratégie optimale $\widehat{\alpha}$   telle que
$$K(\widehat{\alpha},i,x) = \ESS{\alpha \in \A} K(\alpha,i,x).$$
Nous montrons par la suite que notre problème est réduit à prouver l'existence d'un couple de processus $(Y^1,Y^2)$ à l'aide des outils de l'enveloppe de Snell. En effet:
\begin{pro} 
 \label{procY}
Supposons qu'il existe deux processus de  $\mathcal{C}^2$     $Y^1 = (Y_t^1)_{t \geq 0}$ et $Y^2 = (Y_t^2)_{t \geq 0}$ à valeurs dans $\R$  tels que $\forall ~  t \geq 0:$
       \begin{eqnarray}
     Y_t^1 &=& \ESS{\theta \in \mathcal{T}_t}  ~  \mathbb{E} \left [ \int_t^{\theta}e^{- \beta s} f(0,X_s^0) ~ds   -  e^{- \beta \theta} c_{0,1} + Y_{\theta}^2 | \F_t\right],~Y^1_\infty=0   \label{Y1}\\
       Y_t^2 &=& \ESS{\theta \in\mathcal{T}_t}  ~  \mathbb{E} \left [ \int_t^{\theta}e^{- \beta s} f(1,X_s^1) ~ds    - e^{- \beta \theta} c_{1,0}  + Y_{\theta}^1 | \F_t\right],~Y^2_\infty=0\label{Y2}.
      \end{eqnarray}
Alors $Y_0^1 = \sup_{\alpha \in \A} K(\alpha,i,x).$ De plus,  la suite $\widehat{\alpha}= (\widehat{\tau}_n)_{n \geq 0}$  définie par:
       \begin{eqnarray}
   \label{taupair}   \widehat{\tau}_{-1}&=&0\nonumber\\
         \widehat{\tau}_{2n} &=&\inf \{t \geq \widehat{\tau}_{2n-1},  Y_t^1 = - c_{0,1}e^{- \beta t}  +  Y_t^2\}, ~ \forall ~  n \geq 0, \\
       \label{tauimpair}
        \widehat{\tau}_{2n+1} &=&\inf \{t \geq  \widehat{\tau}_{2n}  , Y_t^2 = - c_{1,0}e^{- \beta t}  +  Y_t^1\} 
          \end{eqnarray}
          est optimale  pour le problème de contrôle impulsionnel (\ref{Kalpha}).
 \end{pro}

  \noindent {\bf Preuve. }
  Pour tout $t \geq 0$, nous avons: 
\begin{equation}
\label{1}
    Y_t^1 + \int_0^t e^{- \beta s}  f(0,X_s^0) ~  ds  = \ESS{\theta \in \TT}  ~  \mathbb{E} \left [ \int_0^{\theta}e^{- \beta s} f(0,X_s^0) ~ds   - e^{- \beta \theta}   c_{0,1}  + Y_{\theta}^2 | \F_t\right].
  \end{equation}
Par définition,  $    Y_t^1 + \int_0^t  e^{- \beta s}  f(0,X_s^0)~ds$ est l'enveloppe de Snell du processus 
  \begin{equation}
\label{envS}
U_{\theta}:= \int_0^{\theta}e^{- \beta s} f (\xi_s,X_s^{\xi_s})  ds  - e^{- \beta \theta} c_{0,1}  + Y_{\theta}^2.
 \end{equation}
  Sous la définition  \ref{defopt}, un temps  $\widehat{\tau}_0$ optimal  pour le problème (\ref{envS}) 
   est défini par: 
  $$\widehat{\tau}_0 = \inf \{t \geq  0,     Y_t^1 = -e^{- \beta t}  c_{0,1}  +  Y_t^2 \}.$$
 Par suite,  le temps  $\widehat{\tau}_0$ est     le premier temps défini en (\ref{taupair})  et donc optimal pour le problème (\ref{Y1}) posé. De plus, 
   le processus $Y_0^1$ étant $\F_0$-mesurable, alors $Y_0^1 = \mathbb{E}[Y_0^1].$\\
   Par conséquent, l'égalité (\ref{1}) peut être écrite pour $t=0$ et $\theta = \widehat{\tau}_0$ sous la forme suivante:
   \begin{equation}\label{2}
  Y_0^1   =   \mathbb{E} \left [ \int_0^{\widehat{\tau}_0}  e^{- \beta s} f(0,X_s^0) ~ds   - e^{- \beta \widehat{\tau}_0}   c_{0,1}  + Y_{\widehat{\tau}_0}^2  \right]. 
\end{equation}
En utilisant  la définition de $Y^2$ appliquée à $t = \widehat{\tau}_0$, nous obtenons 
$$
    Y_{\widehat{\tau}_0}^2 = \ESS{\theta \in \mathcal{T},\theta\geq \widehat{\tau}_0}  ~  \mathbb{E} \left [ \int_{\widehat{\tau}_0}^{\theta}e^{- \beta s} f(1,X_s^1) ~ds   - e^{- \beta \theta}   c_{1,0}  + Y_{\theta}^1 | \F_{\widehat{\tau}_0}\right].$$
   De même que     pour  le problème    (\ref{envS}), 
  $  Y_{t}^2 + \int_0^{t}  e^{- \beta s}  f(\xi_s,X_s^\xi) ~  ds $ est l'enveloppe de Snell du processus 
  \begin{equation}
\label{envS1}
U_{\theta}:= \int_0^{\theta}e^{- \beta s}f(\xi_s,X_s^{\xi_s}) ds    - e^{- \beta \theta}   c_{1,0} + Y_{\theta}^1, ~ \theta \geq \widehat{\tau}_0.
 \end{equation}
  Sous la définition  \ref{defopt}, un temps  $\widehat{\tau}_1$ optimal après $\widehat{\tau}_0$  pour le problème (\ref{envS1}) 
   est défini par: 
  $$\widehat{\tau}_1 = \inf \{t \geq \widehat{\tau}_0,     Y_t^2 = -e^{- \beta t}  c_{1,0}  +  Y_t^1 \}.$$
 Par suite,  le temps  $\widehat{\tau}_1$ est     le deuxième temps défini en (\ref{taupair})  et donc optimal pour le problème (\ref{Y2}) posé. D'où, 
 \begin{equation}
\label{3}
  Y_{\widehat{\tau}_0}^2 = \mathbb{E} \left [ \int_{\widehat{\tau}_0}^{\widehat{\tau}_1}e^{- \beta s} f(1,X_s^1) ~ds   - e^{- \beta \widehat{\tau}_1}   c_{1,0}  + Y_{ \widehat{\tau}_1}^1 | \F_{\widehat{\tau}_0}\right],     
\end{equation}
 et  les égalités  (\ref{2}) et (\ref{3}) impliquent: 
  \begin{equation}
   Y_0^1 = \mathbb{E} \left [\int_0^{\widehat{\tau}_0}  e^{- \beta s} f(0,X_s^0) ~ds   +   \int_{\widehat{\tau}_0}^{\widehat{\tau}_1}e^{- \beta s} f(1,X_s^1) ~ds   -
   e^{- \beta\widehat{\tau}_0}  c_{0,1} - e^{- \beta \widehat{\tau}_1}c_{1,0}+ Y_{ \widehat{\tau}_1}^1 \right].\label{5}
     \end{equation}
En utilisant 
$$  \int_{0}^{\widehat{\tau}_1}   e^{- \beta s} f(\xi_s,X_s^{\xi_s}) ~ds =  \int_0^{\widehat{\tau}_0}  e^{- \beta s} f(0,X_s^0) ~ds   +   \int_{\widehat{\tau}_0}^{\widehat{\tau}_1}e^{- \beta s} f(1,X_s^1) ~ds, $$
  l'égalité (\ref{5}) devient:
 $$ Y_0^1 =    \mathbb{E} \left [  \int_{0}^{\widehat{\tau}_1}   e^{- \beta s} f(\xi_s,X_s^{\xi_s}) ~ds  - e^{- \beta\widehat{\tau}_0}  c_{0,1} - e^{- \beta \widehat{\tau}_1}  c_{1,0} +Y_{ \widehat{\tau}_1}^1 \right].$$
  En répétant ce raisonnement successivement, nous obtenons: 
  \begin{equation}\label{6}
  Y_0^1 =    \mathbb{E} \left [  \int_{0}^{\widehat{\tau}_{2 n+1}}   e^{- \beta s} f(\xi_s,X_s^{\xi_s}) ~ds  -  \sum_{k= 0}^n  
  \left( e^{- \beta\widehat{\tau}_{2 k}}  c_{0,1} +  e^{- \beta \widehat{\tau}_{2 k+1}}  c_{1,0} \right)+ Y_{ \widehat{\tau}_{2 n+1}}^1 \right].
  \end{equation}
  Le processus $Y^1$ étant càdlàg et la suite $(\widehat{\tau}_n)$ étant   croissante vers $\widehat{\tau},$ 
   alors les processus $ Y^i_{ \widehat{\tau}_{2 n+1}},~(i=1,2)$ tendent vers $Y^i_{ \widehat{\tau}^-}$ p.s. lorsque $n$ tend vers l'infini. De plus, par définition des processus $Y^1$ et $Y^2,$ nous avons en passant à la limite
  sur l'évènement $\{\tau<\infty\}$:
   $$ Y^i_{ \widehat{\tau}^-} =  - e^{- \beta \widehat{\tau}} c_{0,1}+ Y_{ \widehat{\tau}^-}^{1-i}= - e^{- \beta \widehat{\tau}} (c_{0,1} + c_{ 1,0}) + Y_{ \widehat{\tau}^-}^i.$$
    Il en résulte que $\tau = +\infty$ et puisque  par hypothèse $Y^1_{\infty} = Y^2_{\infty} =0$, en prenant   la limite de  l'égalité (\ref{6}), nous obtenons $Y_0^1 = K(\widehat{\alpha},i,x).$  \\
 
 Montrons maintenant que la stratégie $\widehat{\alpha}$ est optimale, i.e. 
 $$Y_0^1 \geq K(\alpha,i,x), \quad \forall ~\alpha \in \A.$$
   Le temps  $\widehat{\tau}_0$ étant  optimal pour le problème (\ref{Y1}), nous obtenons:
    $$
Y_0^1   \geq   \mathbb{E} \left [ \int_0^ {\tau_0}  e^{- \beta s} f(0,X_s^0) ~ds   - e^{- \beta  \tau_0}  c_{0,1}  + Y_{\tau_0}^2 \right], ~   \forall\tau_0>0.
$$
De plus,   pour tout temps $\tau_1>\tau_0,$ 
$$   Y_{\tau_0}^2 \geq    \mathbb{E} \left [ \int_{\tau_0}^{\tau_1}e^{- \beta s} f(1,X_s^1) ~ds   - e^{- \beta \tau_1}  c_{1,0}  + Y_{\tau_1}^1 | \F_{\tau_0}\right].$$
 Ainsi, nous obtenons pour toute suite croissante   vers l'infini  $(\tau_n)$:
 $$
   Y_0^1 \geq  \mathbb{E} \left [\int_{  0} ^{\tau_1} e^{- \beta s} f(\xi_s,X_s^{\xi_s}) ~ds    -
   e^{- \beta \tau_0}  c_{0,1} - e^{- \beta  \tau_1}  c_{1,0} + e^{- \beta  \tau_1} Y_{  \tau_1}^1 \right]. 
  $$ 
   En répétant successivement le même raisonnement, nous obtenons:
    \begin{equation}\label{7}
     Y_0^1 \geq     \mathbb{E} \left [  \int_{0}^{\tau_{2 n+1}}   e^{- \beta s} f(\xi_s,X_s^{\xi_s}) ~ds  -  \sum_{k= 0}^n 
  \left( e^{- \beta \tau_{2 k}}  c_{0,1} +  e^{- \beta  \tau_{2 k+1}}c_{1,0} \right)+Y_{\tau_{2 n+1}}^1 \right].
  \end{equation}
  La partie droite de l'inégalité (\ref{7}) tend vers $K(\alpha,0,x)+Y_\infty^1=K(\alpha,0,x)$ lorsque n tend vers l'infini.  
  Par conséquent,  $ Y_0^1 = K(\widehat{\alpha},0,x)      \geq K(\alpha,0,x),$ ce qui implique l'optimalité de $\widehat{\alpha}$.
     \endproof

  \section{EDS  rétrogrades et  réfléchies}
  \label{EDSFB}
    Dans cette section, nous étendons au cas de l'horizon infini des  résultats concernant les équations différentielles stochastiques  rétrogrades et réfléchies    sous des hypothèses convenables (cf. \cite{Karoui1} et \cite{hamadene}).\\

\subsection{EDS  rétrogrades}
     
     L'un des résultats fondamentaux concernant les équations différentielles  stochastiques rétrogrades (EDSR) est le théorème donné par  Pardoux et  Peng  \cite{pardoux, pardoux1}, et qui montre l'existence et l'unicité de la solution d'une EDSR en horizon fini sous des hypothèses de Lipschitz sur la fonction drift.   Nous montrons l'existence de la solution en horizon infini en imposant des hypothèses supplémentaires sur la fonction drift et en utilisant des estimations du processus Y.     A la différence des précédents auteurs, notre fonction drift dépend uniquement du temps et du processus $Y.$ 

     \begin{thm} \label{EDSB}
  Supposons   que la fonction $f(.,y)$  est $\F_t$-progressivement mesurable et que:
  $$(\mathcal{H})\left\lbrace
\begin{array}{lcl}
 \forall t, ~y \mapsto f(t,y) ~ \mbox{est décroissante,}\\
 t \mapsto f(t,0)~ \mbox{est bornée,}\\
 \mbox{il existe une constante } C>0 ~ \mbox{telle que  pour tout  } t\geq 0  ~
  \forall y, y'   \in  \R: \\
       |f(t,  y) - f(t,  y') | \leq C   |y-y'|    ~   ~ p.s. 
\end{array}\right.$$            
 Alors, l'équation différentielle stochastique   rétrograde (EDSR)
 \begin{equation}\label{eqB}
 Y_t =  \int_t^{+\infty}  e^{- \beta s}  f(s,Y_s) ds  -  \int_t^{+\infty}  Z_s dW_s, ~ ~Y_{\infty} = 0, ~  t \geq 0
 \end{equation}
 admet une solution $(Y,Z)$ telle que $Y \in \mathcal{C}^2,  ~  Z \in  \mathbb{H}^2$.
   \end{thm}
    
   On montre en préalable un lemme qui sera également utile dans le cas des EDSR réfléchies.
   
    \begin{lem}
   Soit $f$ vérifiant l'hypothèse  $(\mathcal{H})$ et $(Y_t,Z_t,K_t^+, K_t^-),t\in[0,T],$ solution de l'équation différentielle stochastique rétrograde réfléchie   à double barrière:
 {\footnotesize 
 \begin{equation*}
      Y_t =  \int_t^{T}  e^{- \beta s}  f(s,Y_s) ds + \int_t^{T}  e^{- \beta s} dK_s^+ - \int_t^{T}  e^{- \beta s} dK_s^- -  \int_t^{T}  Z_s dW_s , ~   Y_{T} = 0, ~  t \in[0,T],
      \end{equation*}}
où  $Y \in \mathcal{C}^2$,    $Z\in \mathbb{H}^2$, $L_t\leq Y_t \leq   U_t,$
$dK_.^+$ et $dK_.^-$ sont  deux mesures positives vérifiant  $\mathbb{E} [(\int_0^{T}e^{- \beta s} dK_s^+)^2] <\infty,  \mathbb{E} [(\int_0^{T}e^{- \beta s} dK_s^-)^2] <\infty$
     et $$\int_t^{T}   e^{- \beta s}  (Y_s - L_s)dK_s^+  =\int_t^{T}   e^{- \beta s}  (U_s- Y_s )dK_s^-  = 0, ~ ~ \mathbb{P}\mbox{-p.s.}$$
     (i) Sous l'hypothèse $U\geq 0$ et $\sup_tL_t^+\in \mathbb{L}^2,$ et si $f$ ne dépend pas de $y$ on a la majoration
pour tout $\eps>0$
\begin{equation}
  \label{estY_t4}
  \EE[Y_t^2]\mbox{ et }~\EE(\int_t^T|Z_s|^2ds)\leq \varphi(t) \exp\Big(\int_t^T  e^{- \beta u}  du \Big),
   \end{equation}
   où $\varphi(t) =  \frac{1}{\beta}     \|f\|^2  e^{- \beta t} 
  +\frac{1}{\eps}E[sup_{s\geq t}(L_s^+)^2]+  \eps E[(\int_t^T e^{- \beta s}dK_s^+)^2].$\\
\\
(ii) Si $f$ dépend de $y$ avec $C$  la constante de Lipschitz de la fonction $y\mapsto f(s,y)$, mais si $L\leq 0 \leq U$ voire si $L=U=0$ (c'est à dire le cas non réfléchi), alors pour tout  $ t \in[0,T]:$  
    \begin{equation}
  \label{estY_t3}
   \EE[Y_t^2]\leq 
  D e^{- \beta t},
   \end{equation}
   où 
$D =  \frac{1}{\beta} \| f\|^2  \exp\Big( \frac{2 C+1}{\beta }  \Big). $
   \end{lem}

  \noindent{\textbf{Preuve.}} C'est une adaptation de la preuve de la Proposition 3.5 de \cite{Karoui1}: 
  la formule d'Itô
  et $\int_t^{T}   e^{- \beta s}  (Y_s - L_s)dK_s^+  =\int_t^{T}   e^{- \beta s}  (U_s- Y_s  )dK_s^-  = 0$ montrent:
  {\footnotesize
  $$Y_t^2+\int_t^T|Z_s|^2ds= 2\int_t^{T} Y_s e^{- \beta s}  f(s,Y_s) ds + 2\int_t^{T} L_s e^{- \beta s} dK_s^+ - 2\int_t^{T} U_s e^{- \beta s} dK_s^- - 2 \int_t^{T}  Y_sZ_s dW_s .
  $$}
  (i) Dans ce premier cas, $f$ ne dépend pas de $y$ et de ce fait est bornée par $\|f\|$ 
  et 
  $$2\int_t^{T} Y_s e^{- \beta s}  f(s,\omega) ds\leq \int_t^{T} Y_s^2 e^{- \beta s} ds+\int_t^{T}  e^{- \beta s}  \|f\|^2 ds,$$
   de plus,   pour tout $\eps>0:$
\begin{eqnarray*}
2\int_t^{T} L_s e^{- \beta s} dK_s^+ &\leq& 2\int_t^{T} L_s^+ e^{- \beta s} dK_s^+\nonumber\\
&\leq&
   2\sup_{s\geq t}L_s^+\int_t^{T} e^{- \beta s} dK_s^+\leq 
   \frac{1}{\eps}sup_{s\geq t}(L_s^+)^2+
  \eps\Big(\int_t^T e^{- \beta s}dK_s^+\Big)^2
  \end{eqnarray*}
   (on utilise $2ab\leq \frac{1}{\eps}a^2+\eps b^2$).
  Ainsi, puisque $U\geq 0$ et $dK^-$ mesure positive, il vient globalement
  {\footnotesize
  $$\EE[Y_t^2]+\int_t^T|Z_s|^2ds]\leq \int_t^{T}\EE[Y_s^2] e^{- \beta s} ds+\frac{1}{\beta}   e^{- \beta t}  \|f\|^2 ds
  +\frac{1}{\eps}\EE[sup_{s\geq t}(L_s^+)^2]+  \eps \EE\Big[\Big(\int_t^T e^{- \beta s}dK_s^+\Big)^2\Big] .$$}
  On applique le lemme de Gronwall  pour majorer la fonction
  $t\mapsto \EE[Y_t^2]$ (cf. appendice) avec  
  $\psi (t) =  e^{- \beta t}$ et 
\begin{equation}
\label{def:varphi}
\varphi(t) =  \frac{1}{\beta}     \|f\|^2  e^{- \beta t} 
  +\frac{1}{\eps}\EE[sup_{s\geq t}(L_s^+)^2]+  \eps \EE\Big[\Big(\int_t^T e^{- \beta s}dK_s^+\Big)^2\Big].
  \end{equation}
   Comme   $\varphi$ est décroissante il vient:
  $$\EE[Y_t^2]\leq \varphi(t)\exp\Big(\int_t^T  e^{- \beta u}  du \Big).$$
  On a aussi par conséquent, puisque $\int_t^T|Z_s|^2ds\leq Y_t^2+\int_t^T|Z_s|^2ds,$
  la même majoration
  $$\EE(\int_t^T|Z_s|^2ds)\leq \varphi(t)\exp\Big(\int_t^T  e^{- \beta u}  du \Big).$$
  \\
(ii) Dans le cas où $f$ est Lipschitzienne de coefficient $C$
  et $L\leq 0\leq U$ on reprend le développement de Ito:
  {\footnotesize
  $$Y_t^2+\int_t^T|Z_s|^2ds= 2\int_t^{T} Y_s e^{- \beta s}  f(s,Y_s) ds + 2\int_t^{T} L_s e^{- \beta s} dK_s^+ - 2\int_t^{T} U_s e^{- \beta s} dK_s^- - 2 \int_t^{T}  Y_sZ_s dW_s .
  $$}
  dont on prend l'espérance après avoir utilisé la propriété de Lipschitz de la fonction $y\mapsto f(s,y)$, la majoration $2| Y_s||  f(s,0)|\leq | Y_s|^2+|  f(s,0)|^2$ et   $L_t \leq 0 \leq U_t:$
  \begin{eqnarray*}
  E[Y_t^2]&\leq&\EE \left[ 2C\int_t^{T} Y_s^2 e^{- \beta s} ds +\int_t^{T}e^{- \beta s} 2| Y_s||  f(s,0)| ds  \right]\nonumber\\
  &\leq& \EE \Big[(2C+1)\int_t^{T} Y_s^2 e^{- \beta s} ds+\int_t^{T}e^{- \beta s}  f^2(s,0) ds \Big].
  \end{eqnarray*}
    On applique le lemme de Gronwall  (cf. appendice) avec   $$\varphi (t) = \frac{1}{\beta} \| f\|^2  e^{- \beta t} ~ \mbox{ et } ~ ~ 
 \psi (t) = (2C+1) e^{- \beta t}.$$
  Comme  $\varphi$ est décroissante,  il vient:
 
\begin{eqnarray*}
 \EE[Y_t^2]  &\leq& \frac{1}{\beta} \| f\|^2  e^{- \beta t} 
 \exp\Big(  (2 C+1)   \int_t^T   e^{- \beta u}du \Big) \nonumber\\
 &\leq&  \frac{1}{\beta} \| f\|^2  e^{- \beta t}  \exp\Big( \frac{2 C+1}{\beta } \Big ).
 \end{eqnarray*}

\endproof

  \noindent{\textbf{Preuve du théorème \ref{EDSB}. }}
     En utilisant la proposition 2.2 de \cite[p. 57]{pardoux}, 
    on déduit  l'existence 
    $\forall n$ d'un couple de processus $(Y^n,Z^n)$ qui vérifient $Y^n \in \mathcal{C}^2,~  Z^n \in  \mathbb{H}^2$  et 
    $ ~ \forall t \leq n:$  
   $$Y^n_t =  \int_t^n e^{- \beta s}  f(s,Y_s^n) ds  -  \int_t^{n}  Z_s^n dW_s, ~ Y_t^n=0  ~ \forall ~  t \geq n.$$
  Ensuite, $\forall k  \geq n,$ nous obtenons:
  
  \begin{eqnarray*}
 Y^{n+k}_t - Y^n_t   &=&  \int_t^n e^{- \beta s} \left[f(s,Y^{n+k}_s) -  f(s,Y_s^n)  \right  ] ds  + \int_n^{n+k} e^{- \beta s}  f(s,Y^{n+k}_s)    ds 
 \nonumber\\
 & - &  \int_t^{n}\left[  Z_s^{n+k} - Z_s^{n} \right  ]  dW_s - \int_n^{n+k}  Z_s^{n+k} dW_s.
     \end{eqnarray*}
  Appliquons la formule d'Itô à $( Y^{n+k}_t-Y^n_t)^2$ entre $t$ et $n$, 
   où l'on ajoute le développement de $Y^{n+k}$ entre $n$ et $n+k$ 
  nous obtenons:
   
  $$ ( Y^{n+k}_t - Y^n_t  )^2  + \int_t^{n}\left(  Z_s^{n+k} - Z_s^{n} \right  )^2  ds + \int_n^{n+k}  (Z_s^{n+k}) ^2 ds$$
{\footnotesize   \begin{eqnarray*}
  &=& 2 \int_t^n e^{- \beta s} ( Y^{n+k}_s - Y^n_s  ) \left[f(s,Y^{n+k}_s) -  f(s,Y_s^n)  \right  ] ds  + 2\int_n^{n+k} e^{- \beta s}  Y^{n+k}_s   f(s,Y^{n+k}_s)    ds 
 \nonumber\\
 & - &  2\int_t^{n}( Y^{n+k}_s - Y^n_s  ) \left(  Z_s^{n+k} - Z_s^{n} \right  )  dW_s - 2\int_n^{n+k}   Y^{n+k}_s    Z_s^{n+k} dW_s.
     \end{eqnarray*}}
  En utilisant que $\forall y, ~ t \mapsto f(t,y) ~ \mbox{est décroissante}$ et en passant à l'espérance, nous obtenons:
  {\footnotesize
  $$ \mathbb{E}( Y^{n+k}_t - Y^n_t  )^2  + \mathbb{E} \int_t^{n}\left(  Z_s^{n+k} - Z_s^{n} \right  )^2  ds + \mathbb{E}\int_n^{n+k}   (Z_s^{n+k}) ^2 ds \leq 
   2~ \mathbb{E} \int_n^{n+k} e^{- \beta s}   Y^{n+k}_s     f(s,Y^{n+k}_s) ds
  $$
  }
\begin{eqnarray*}
&\leq &  2~ \mathbb{E} \int_n^{n+k} e^{- \beta s} | Y^{n+k}_s  | |f(s,Y^{n+k}_s) - f(s,0)|ds +  2\mathbb{E} \int_n^{n+k} e^{- \beta s} | Y^{n+k}_s  | |  f(s,0)|ds. 
    \end{eqnarray*}
    En utilisant la condition de Lipschitz, il vient: 
    {\footnotesize  \begin{eqnarray*}
     \mathbb{E}( Y^{n+k}_t - Y^n_t  )^2  &+& \mathbb{E} \int_t^{n}\left(  Z_s^{n+k} - Z_s^{n} \right  )^2  ds 
     \\
    &\leq &  2\mathbb{E} \int_n^{n+k}   C e^{- \beta s} | Y^{n+k}_s  |^2 ds +  2\mathbb{E} \int_n^{n+k}e^{- \beta s} | Y^{n+k}_s  | |  f(s,0)|ds\nonumber
      \\
      &\leq &
      2\mathbb{E} \int_n^{n+k}   C e^{- \beta s} | Y^{n+k}_s  |^2 ds +  2\left[ \mathbb{E}  \int_n^{n+k} e^{- \beta s} | Y^{n+k}_s  |^2\right]^\demi
      \left[ \mathbb{E}  \int_n^{n+k} e^{- \beta s} |  f(s,0)|^2 ds \right]^\demi,
     \end{eqnarray*}}
   où la dernière égalité est obtenue en utilisant l'inégalité de Cauchy-Schwartz. De plus,  $s\mapsto f(s,0)$ étant bornée, nous avons: 
    $$\left[ \mathbb{E}  \int_n^{n+k} e^{- \beta s} |  f(s,0)|^2 ds \right]^\demi \To 0$$
     lorsque n tend vers l'infini. 
     Ensuite, en utilisant (\ref{estY_t3}) pour tout $n$ et tout $k$ avec $L=U= 0:$
  $$
   \EE[(Y_s^{n+k})^2]\leq D e^{- \beta s}, 
$$
 où $D = \frac{1}{\beta} \| f\|^2   \exp\Big( \frac{2 C+1}{\beta } \Big ).$
    En passant à l'intégrale et en utilisant Tonelli et la bornitude de $s\mapsto f(s,0)$, 
    $$ \mathbb{E} \int_n^{n+k} e^{- \beta s} | Y^{n+k}_s  |^2 ds  \leq   D \int_n^{n+k} e^{- 2\beta s}ds 
  \leq 
\frac{D}{2 \beta} e^{-2 \beta n}. 
    $$
   Ainsi,  lorsque $n$ tend vers l'infini
    $$\mathbb{E} \int_n^{n+k} e^{- \beta s} | Y^{n+k}_s  |^2 ds \To 0  $$
    lorsque n tend vers l'infini. Par conséquent,  les suites  $(Y^n)$ et $(Z^n)$ sont deux suites de Cauchy qui convergent   respectivement dans $\mathbb{L}^2 ( \Omega)$
    et $\mathbb{L}^2 ([0,\infty[ \times \Omega, dt \otimes d\mathbb{P})$ vers deux processus $Y, Z$
    qui sont donc respectivement élément de $\mathbb{L}^2 ( \Omega)$
    et $\mathbb{L}^2 ([0,\infty[ \times \Omega, dt \otimes d\mathbb{P})$: 
    $$\lim_{n \To \infty}Y^n_t = Y_t ~ ~ \mbox{et} ~ ~ \lim_{n \To \infty}Z^n_t = Z_t.$$
   On a donc  pour tout $t, ~ Y_t \in \mathbb{L}^2 ( \Omega),$
   $  Z \in \mathbb{L}^2 ([0,\infty[ \times \Omega, dt \otimes d\mathbb{P}),~ \EE (Y_t)^2 = \displaystyle\lim_{n \To \infty } \EE (Y_t^n)^2$
 et par suite: 
 $$
 \EE (Y_t)^2 \leq D e^{-\beta t}.
  $$
  
Par limite presque sûre     de la relation  vérifiée pour tout $n$ par le couple
   $(Y^n,Z^n),$ plus la continuité sur $\R$  de la fonction $f$,   on obtient L'EDSR (\ref{eqB}).\\

  Enfin, de la relation (\ref{eqB})
    on tire 
$$      Y_t = Y_0- \int_0^t e^{- \beta s}  (f(s,Y_s)-f(s,0)) ds - \int_0^t e^{- \beta s}  f(s,0) ds  +  \int_0^t Z_s dW_s ,$$
le deuxième terme  avec la condition de Lipschitz est majoré en valeur absolue par
$$\big|\int_0^t e^{- \beta s}  (f(s,Y_s-f(s,0)) ds \big|\leq C\int_0^t e^{- \beta s}  |Y_s|ds,$$ et
 $$\sup_{  t}|Y_t|\leq |Y_0|+C\int_0^{\infty} e^{- \beta s}  |Y_s|ds+\frac{\|f\|}{\beta} + \sup_{  t}| \int_0^{\infty} Z_s dW_s|,$$
c'est à dire   une majoration par la somme de quatre termes dans $\LL^2:$
$\EE( \sup_{  t}|Y_t|^2)   <\infty,$ et  $Y \in \mathcal{C}^2.$ 

\endproof

\begin{cor} \label{compsp}
Supposons que   $(Y,Z)$ et $(Y',Z')$  sont  solutions de l'EDS rétrograde (\ref{eqB}) associées à $f$ et $f'$, avec $f$   et $f'$ satisfaisant $(\mathcal{H})$. 
Supposons de plus que $ \forall y \in \R$:  
$$(\mathcal{H}_1): \quad    f(t,y) \leq f'(t,y)~  ~ dt \otimes d\mathbb{P} \mbox{-p.s.} $$ Alors: 
  $$Y_t \leq Y_t', ~    ~ t\geq 0, ~  ~ \mathbb{P}\mbox{-p.s.}$$
     \end{cor}
     \noindent{\textbf{Preuve.}} La preuve de ce corollaire est une conséquence immédiate de l'hypothèse $(\mathcal{H}_1)$  et du  théorème 4.1 de    El Karoui et al.  \cite[p. 712]{Karoui1}. 
     En effet, sous $(\mathcal{H}_1)$, nous obtenons:
      $$ Y_t^n \leq Y_t'^n, ~    ~\forall ~ n,  ~ \forall t\in[ 0,n],  ~  ~ \mathbb{P} \mbox{-p.s.}$$     
     En faisant tendre $n$ vers l'infini nous obtenons la comparaison recherchée: $$Y_t \leq Y_t', ~    ~ t\geq 0, ~  ~ \mathbb{P}\mbox{-p.s.}$$
     \endproof 
    
       \subsection{EDS  rétrogrades  réfléchies }
        Dans cette section, nous généralisons des résultats de  El Karoui et al. \cite{Karoui1} et  Hamadène et al. \cite{hamadene}  pour la résolution des  équations différentielles stochastiques  rétrogrades réfléchies (EDSRR) à double barrière. 
 On   dispose de deux données:\\
   \\
  $\bullet$ Rappelons que la fonction $f(.,y)$  est $\F_t$-progressivement mesurable et que:
  $$(\mathcal{H})\left\lbrace
\begin{array}{lcl}
 \forall t, ~y \mapsto f(t,y) ~ \mbox{est décroissante,}\\
 t \mapsto f(t,0)~ \mbox{est bornée,}\\
 \mbox{il existe une constante } C>0 ~ \mbox{telle que pour tout } ~ t \geq 0, ~ 
  \forall y, y'   \in  \R: \\
       |f(t,  y) - f(t,  y') | \leq C   |y-y'|    ~   ~ p.s. 
\end{array}\right.$$            
\\
      $\bullet$    Une  barrière $(L_t)_{t\geq 0}$: un processus
          continu     à valeurs réelles $\F_.$-adapté satisfaisant:
         $$  \EE \Big ( \sup_{t \geq 0}  (L_t^+)^2 \Big) < \infty.     $$
        
      Nous étendons au cas de l'horizon infini    la proposition 5.1  de  El Karoui et al. \cite[p. 716]{Karoui1}: nous  étudions des EDS rétrogrades  réfléchies   en horizon infini avec  une  barrière $(L_t)_{t\geq 0}$.  
    \begin{thm}
    \label{bar1}
         Sous les hypothèses  $(\mathcal{H}) $    mais avec $f$ ne dépendant pas
         de $y$  et    $\EE \Big ( \sup_{t \geq 0}  (L_t^+)^2 \Big) < \infty ,$   il existe une solution $(Y,Z,K)$ telle que pour tout $t \geq 0$: 
      \begin{itemize}
    \item[i/] $Y \in \mathcal{C}^2$, $Z\in \mathbb{H}^2$. 
      \item[ii/]  \begin{equation}
      \label{EDSR1}
      Y_t =  \int_t^{+\infty}  e^{- \beta s}  f(s) ds + \int_t^{+\infty}  e^{- \beta s} dK_s -  \int_t^{+\infty}  Z_s dW_s , ~   Y_{\infty} = 0. 
      \end{equation}
 \item[iii/] $Y_t \geq   L_t.$
   \item[iv/] $(dK_t)$ est   une mesure positive vérifiant   $\mathbb{E} [(\int_0^\infty e^{- \beta s} dK_s)^2] <\infty$ 
               et $$\int_0^t   e^{- \beta s}  (Y_s - L_s)dK_s  = 0, ~ ~ \mathbb{P}\mbox{-p.s.}$$
   \item[v/]  
                 $Y_t = \ESS{\theta\geq t } ~  \mathbb{E} \Big[\int_t^\theta e^{- \beta s} f(s,\omega)ds + L_{\theta}|\F_t \Big].$
 \end{itemize}
 \end{thm}

  \noindent{\textbf{Preuve. }} $i/$ En utilisant la proposition 5.1 de \cite[p. 716]{Karoui1},
   on déduit
      l'existence 
    $\forall n$ d'un triplet  de processus $(Y^n,Z^n,K^n)$ qui vérifient $Y^n \in \mathcal{C}^2([0,n]),~  Z^n \in  \mathbb{H}^2 ([0,n] \times \Omega)$  et 
    $ ~ \forall t \leq n:$  
         
         $$ \left\lbrace
\begin{array}{lcl}
 &&Y^n_t \geq L_t,~ Y_t^n=0~~\forall t\geq n. \\
 \\
 && Y^n_t   =   \int_t^n e^{- \beta s}  f(s) ds +  \int_t^n e^{- \beta s}  dK^n_s -  \int_t^{n}  Z_s^n dW_s. \\
  \\
 &&\mathbb{E} (\int_0^n e^{- \beta s}  dK_s^n )^2 < \infty ~ \mbox{et} ~ \int_0^t   e^{- \beta s}  (Y_s^n - L_s)dK_s^n  = 0, ~ ~  t \leq n  ~ ~  \mathbb{P}\mbox{-p.s.}
\end{array}\right.$$    
     
    Le théorème de comparaison de S. Hamadène et al.  (\cite[p. 164]{hamadene}, proposition 41.3) appliquée à la suite croissante   $f_n(s) = f(s) \1_{[0,n]}$ implique que $(Y^n)_n$ est une suite croissante de processus et que $(\int e^{- \beta s} dK^n )_n$ est une suite décroissante de processus. On note leurs limites
   presque sûres  respectivement  $Y$ et $\int e^{- \beta s} dK_s.$  \\
 \\
  Ensuite, $\forall k  \geq 0$ et $\forall t\leq n$ nous obtenons:
  {\footnotesize
  \begin{eqnarray}
    \label{SDC} 
 &&Y^{n+k}_t - Y^n_t   =  \int_n^{n+k} e^{- \beta s}  f(s,\omega)ds +   \int_t^n e^{- \beta s} \left[ dK^{n+k}_s  -dK^{n}_s  \right  ] 
\\
 &&   + \int_n^{n+k} e^{- \beta s}dK^{n+k}_s   -  \int_t^{n}\left[  Z_s^{n+k} - Z_s^{n} \right]  dW_s - \int_n^{n+k}  Z_s^{n+k} dW_s.
  \nonumber
     \end{eqnarray}
     }
  Appliquons la formule d'Itô à $( Y^{n+k}_t - Y^n_t  )^2$ entre $t$ et $n$,  nous obtenons:
   \begin{eqnarray*}
 && ( Y^{n+k}_t - Y^n_t  )^2  + \int_t^{n}\left(  Z_s^{n+k} - Z_s^{n} \right  )^2  ds =(Y_n^{n+k})^2
  \\
 &+&   2\int_t^{n}e^{- \beta s}( Y^{n+k}_s - Y^n_s  ) \left(  dK_s^{n+k} - dK_s^{n} \right  )  
  -   2\int_t^{n}( Y^{n+k}_s - Y^n_s  ) \left(  Z_s^{n+k} - Z_s^{n} \right  )  dW_s.
     \end{eqnarray*}
 En utilisant    la croissance de la suite $(Y^n_s)$
  et la décroissance de   la suite $(\int e^{- \beta s} dK^n )_n$, puis en prenant l'espérance, nous obtenons:
  $$
   \mathbb{E}( Y^{n+k}_t - Y^n_t  )^2 
 +\mathbb{E} \int_t^{n}\left(  Z_s^{n+k} - Z_s^{n} \right  )^2  ds
 \leq \EE(Y_n^{n+k})^2.
 $$
 En utilisant (\ref{estY_t4}), nous obtenons pour tout  $\eps>0,$ $n$ et tout  $k$: 
 
 \begin{equation} \label{estY_t5} 
\EE(Y^{n+k}_n)^2\mbox{ et }\EE[(\int_n^{n+k} Z^{n+k}_s dW_s)^2]  \leq   
\varphi(n)  \exp\Big (\int_n^{n+k}  e^{- \beta u}  du\Big)  \leq \exp \Big( \frac{1}{\beta}\Big)\varphi(n)
   \end{equation}
 rappelant (\ref{def:varphi})
   $$\varphi(n)= \frac{1}{\beta}     \|f\|^2  e^{- \beta n} 
  +\frac{1}{\eps}E[sup_{s\geq n}(L_s^+)^2]+  \eps E[(\int_n^{n+k}e^{- \beta s}dK_s^{n+k})^2].$$
    De l'équation  
  $$ \int_n^{n+k}  e^{- \beta s } dK^{n+k}_s  = Y^{n+k}_n - \int_n^{n+k}  e^{- \beta s } f(s,\omega) ds + \int_n^{n+k} Z^{n+k}_s dW_s,$$
   et de l'estimation (\ref{estY_t5}), nous obtenons:
    \begin{equation*} 
\frac{1}{3}\EE\big(\int_n^{n+k}  e^{- \beta s } dK^{n+k}_s \big)^2  \leq  
2\exp \Big( \frac{1}{\beta}\Big) \varphi(n)+e^{- \beta n }\frac{1}{\beta}     \|f\|^2.
    \end{equation*}  
    Si on retranche de $\varphi(n)$  le terme $\eps \EE[(\int_n^{n+k}e^{- \beta s}dK_s^{n+k})^2]$
    il vient
  {\footnotesize    \begin{equation}
   \label{majIntK}
   \big(\frac{1}{3}-2\eps\exp \Big( \frac{1}{\beta}\Big) \big)\EE\Big(\int_n^{n+k}  e^{- \beta s } dK^{n+k}_s \Big)^2  \leq  
    \big(1 +2  \exp \Big( \frac{1}{\beta}\Big)\big )\frac{1}{\beta} \|f\|^2  e^{- \beta n} 
  +\frac{1}{\eps}\EE[\sup_{s\geq n}(L_s^+)^2], 
   \end{equation}}
 donc tend vers zero uniformément dès que $\eps$ est choisi assez petit: en effet, 
   puisque  $\sup_{s}L_s^+\in \LL^2,$  par le théorème de Lebesgue de convergence monotone
    $\EE[sup_{s\geq n}(L_s^+)^2] $ tend vers $0$ quand $n$ tend vers l'infini.
    Globalement $\varphi(n) \rightarrow 0$ quand $n$ tend vers l'infini et on obtient par (\ref{estY_t5}) 
     que     les suites  $(Y^n)$ et $(Z^n)$ sont deux suites de Cauchy qui convergent
       respectivement dans $\mathbb{L}^2(\Omega)$ et $\mathbb{L}^2 ([0,\infty[ \times \Omega, dt \otimes d\mathbb{P})$ vers deux processus $Y, Z:$  
    $$\lim_{n \To \infty}Y^n_t = Y_t ~ ~ \mbox{et} ~ ~ \lim_{n \To \infty}Z^n_t = Z_t.$$
    On a donc $Y_t\in \mathbb{L}^2$ pour tout $t$, $Z\in \mathbb{L}^2 ([0,\infty[ \times \Omega, dt \otimes d\mathbb{P})$ et 
    $\EE(Y_t^2) = \displaystyle\lim_{n\To \infty}\EE(Y_t^n)^2,$ d'où, 
    \begin{equation}\label{majorD''}
    \EE(Y_t^2) \leq \exp \Big( \frac{1}{\beta}\Big) \varphi(t).
    \end{equation}
  $ii/ $  Par limite presque sûre     de la relation  vérifiée pour tout $n$ par le triplet 
   $(Y^n,Z^n, K^n),$     on obtient L'EDSRR (\ref{EDSR1}).\\
\\
   $iii/$ Puisque pour tout $n,$ $Y^n_t\geq L_t$ à la limite on a également $Y_t\geq L_t$
   presque sûrement soit l'item $iii/$.\\
    \\
  %
$iv/ $ Examinons maintenant la différence pour    tout $n$ et $0\leq t\leq n\leq n+k$
  $$\int_t^ne^{-\beta s}dK^n_s-\int_t^{n+k}e^{-\beta s}dK^{n+k}_s=
  (\int_t^ne^{-\beta s}dK^n_s-\int_t^{n}e^{-\beta s}dK^{n+k}_s)-\int_n^{n+k}e^{-\beta s}dK^{n+k}_s$$
  dont le second terme tend vers $0$ dans $\LL^2$ par (\ref{majIntK}).
  Quant au premier utilisant (\ref{SDC}) il se récrit
  $$\int_t^ne^{-\beta s}dK^n_s-\int_t^{n}e^{-\beta s}dK^{n+k}_s=(Y^n_t-Y_t^{n+k})+Y_n^{n+k}+\int_t^n(Z_s^n-Z_s^{n+k})dW_s.$$
  Ces trois termes convergent vers $0$ dans $\LL^2$ : 
  le premier et le troisième sont les restes de Cauchy des suites de Cauchy définissant respectivement
  $Y$ et $Z$. Le second vérifie (\ref{estY_t5}):
  $\EE[(Y_n^{n+k})^2]\leq \exp \big(\frac{1}{\beta }\big)\varphi(n).$
Donc   $(\int_t^ne^{-\beta s}dK^n_s)_n$ est une suite de Cauchy dans $\LL^2$, et sa limite
$\int_t^\infty e^{-\beta s}dK_s$ est donc élément de $L^2.$\\
  \\
  Pour achever la preuve de   $iv/$, pour tout $n$ on a l'identité
 $$\int_0^t  e^{-\beta s}(Y_s^n-L_s)dK_s^n=0,  ~  ~ t \leq n.$$
    La monotonie de la suite de mesures $K^n$ montre que $\lim_n~p.s.\int_0^t e^{-\beta s}L_sdK_s^n=
\int_0^t e^{-\beta s}L_sdK_s. $
    Par ailleurs développons  (en rappelant que par définition $Y_s^n=0$ pour tout $s\geq n$) la différence suivante:
    $$\int_0^te^{-\beta s}Y_s^ndK_s^n-\int_0^t e^{-\beta s}Y_sdK_s=$$
    $$\int_0^te^{-\beta s}Y_s^n(dK_s^n-dK_s)+\int_0^t e^{-\beta s}(Y^n_s-Y_s)dK_s.$$
    Le deuxième terme tend presque sûrement vers $0$ par convergence de Lebesgue monotone.
    Pour le premier, notons que $L_s\leq Y_s^n\leq Y_s$, il converge vers 0 par convergence de Lebesgue
    majorée.\\
    \\
     Suite de l'item $i/$: De la relation (\ref{EDSR1})
    on tire 
$$      Y_t = Y_0- \int_0^t e^{- \beta s}  f(s) ds -\int_0^t  e^{- \beta s} dK_s +  \int_0^t Z_s dW_s.$$
on a donc puisqu'ici  $f$ ne dépend pas de $y$ et dans ce cas  est bornée :
$$\sup_t|Y_t|\leq |Y_0|+\|f\|\frac{1}{ \beta }   +\int_0^\infty  e^{- \beta s} dK_s + \sup_t| \int_0^\infty  Z_s dW_s| $$
c'est à dire une majoration par la somme de quatre termes dans $\mathbb{L}^2$ : 
 $\EE[\sup_t|Y_t|^2]<\infty$ et $Y\in \mathcal{C}^2.$\\
\\
  $v/$   En utilisant la preuve de la proposition 5.1 de \cite[p. 716]{Karoui1}, l'hypothèse
  $\sup_tL_t^+\in \mathbb{L}^2$ étant vérifiée, on a de plus pour tout $n$: 
    $$Y_t^n = \ESS{t \leq \theta \leq n} ~  \mathbb{E} \Big[\int_t^\theta e^{- \beta s} f(s)ds + L_{\theta}|\F_t \Big].$$
   On obtient donc l'item $v/$ par simple limite croissante presque sûre des deux membres de l'égalité.
%
%
   \endproof

 \begin{cor}
 \label{comp}  
   Supposons que    $(Y,Z,K)$ et $(Y',Z',K')$  sont solutions de l'équation stochastique   rétrograde  réfléchie (\ref{EDSR1}) associées à $(f,L)$ et $(f',L)$,  satisfaisant $(\mathcal{H})$ et  $(\mathcal{H}_1)$. Alors:
   $$\forall t\geq 0,~Y_t \leq Y_t' ~ \mbox{ et } ~  \int_0^t e^{- \beta s} dK_s\geq \int_0^t e^{- \beta s} dK_s',~ ~ \mathbb{P}\mbox{-p.s.}$$
 \end{cor}
 \noindent{\textbf{Preuve. }}  
  Soient  $(Y^n,Z^n) $  et $(Y^{'n},Z^{'n}) $    solutions des   EDS rétrogrades ordinaires 
  suivantes: 
    \begin{eqnarray*} 
    Y^n_t & = & \int_t^T  e^{- \beta s}  f(s) ds+ n \int_t^T  e^{- \beta s}  (Y^n_s - L_s)^- ds -  \int_t^T  Z_s^n dW_s, ~ Y_T^n =0  \nonumber\\
    Y^{'n}_t & = & \int_t^T  e^{- \beta s}  f'(s) ds+ n \int_t^T  e^{- \beta s}  (Y^{'n}_s - L_s)^- ds -  \int_t^T  Z_s^{'n} dW_s, ~ Y_T^{'n} =0.
   \end{eqnarray*}
   On pose $f_n(t,y) = f(t) - n   e^{- \beta s}  (y - L_s)^- $ et $f_n'(t,y) = f(t) - n   e^{- \beta s}  (y - L_s)^- $. 
  Sous $(\mathcal{H}_1)$, nous obtenons   $f_n \leq f_n'$ et par suite, en utilisant le résultat  de comparaison des EDS rétrogrades (corollaire \ref{compsp} ), on a 
  $ Y_t^n \leq Y_t'^n, ~ ~\forall ~ n,  t\geq 0, ~  ~ \mathbb{P} \mbox{-p.s.} $
  Ensuite, en utilisant   El Karoui et al. \cite[p. 719]{Karoui1}, on a: 
  $$     \int_t^T  e^{- \beta s} dK^n_s = n \int_t^T  e^{- \beta s}  (Y^n_s - L_s)^- ds, $$
  et en utilisant $ Y_t^n \leq Y_t'^n$, il vient: 
  $$n \int_t^T  e^{- \beta s}  (Y^n_s - L_s)^- ds \geq n \int_t^T  e^{- \beta s}  (Y^{'n}_s - L_s)^- ds.$$
  Par conséquent:
      $$  
      \int_0^t  e^{- \beta s} dK_s^n\geq \int_0^t  e^{- \beta s} dK_s^{'n},~\forall t\leq n, ~  ~ \mathbb{P} \mbox{-p.s.   }$$ 
     En faisant tendre $n$ vers l'infini nous obtenons les comparaisons recherchées:
      $$Y_t \leq Y_t', ~    ~ t\geq 0, ~  ~ \mathbb{P}\mbox{-p.s. 
     et }\int_0^t e^{- \beta s} dK_s\geq \int_0^t e^{- \beta s} dK_s',~ ~ \mathbb{P}\mbox{-p.s.}$$
  \endproof

      \subsection{EDSR réfléchies à double barrière }
      
 Dans cette section, nous généralisons le théorème 42.2 de  Hamadène et al. \cite[p. 167]{hamadene} dans le cas de l'horizon infini.
 Supposons maintenant que notre système admet deux    barrières 
 $(L_t)_{t\geq 0}$ et $(U_t)_{t\geq 0}$: deux processus
           continus    à valeurs réelles $\F_.$-adaptés satisfaisant:
        $$  L_t\leq 0 \leq U_t.$$

      \begin{thm}\label{bar2}
      Sous les hypothèses     $({\cal H})$ et  $L_t\leq 0 \leq U_t,$ 
       il existe un processus  $(Y,Z,K^+,K^-)$ tel que pour tout $t \geq 0$:
      \begin{itemize}
    \item[i/]   $ Y \in   \mathcal{C}^2$  et   $Z\in \mathbb{H}^2$. 
      \item[ii/]  \begin{equation}\label{EDSR2}
      Y_t =   \int_t^{+\infty}  e^{- \beta s}  f(s) ds  + \int_t^{+\infty}  e^{- \beta s} dK_s^+ - \int_t^{+\infty}  e^{- \beta s} dK_s^-  -  \int_t^{+\infty}  Z_s dW_s, ~ Y_{\infty}=0.
      \end{equation}
              \item[iii/] $  L_t \leq Y_t \leq     U_t.$
              \item[iv/] $(dK^+_t)$ et $(dK^-_t)$ sont deux mesures positives  vérifiant   $\mathbb{E} (\int_0^\infty e^{- \beta s}  dK_s^{+} )^2 < \infty,$ $ \mathbb{E} (\int_0^\infty e^{- \beta s}  dK_s^{-} )^2 < \infty$ et 
              $$\int_0^t (Y_s -     L_s)e^{- \beta s}dK_s^+  = \int_0^t
               ( U_s -Y_s )e^{- \beta s}dK_s^-= 0, ~ ~ \mathbb{P}\mbox{-p.s.}$$ 
                                  \end{itemize}
                                  
 \end{thm}
 
 On montre en préalable des  résultats  qui seront  utiles  pour démontrer le théorème \ref{bar2}. 
En utilisant \cite{hamadene,martinLepl}, on déduit   une comparaison des processus dans le cas des équations différentielles réfléchies à double barrière en horizon fini:

\begin{pro}\label{comp2bar}
 Soient   $(Y_t,Z_t,K^+_t,K^-_t)$ et $(Y_t',Z_t',K_t'^+,K_t'^-), ~ t \in[0,T]    $     solutions des EDS rétrogrades  réfléchies  associées à $(f,L,U)$ et $(f',L,U)$  où $f$ ne dépend pas de $y$  satisfaisant $(\mathcal{H})$ et $(\mathcal{H}_1)$. Alors:
 $$\forall t  \in[0,T] , ~  Y_t \leq Y_t',  ~   ~  \int_0^t e^{- \beta s}dK_s^+ \geq \int_0^t e^{- \beta s}
 dK'^+ _s ~ \mbox{et} ~ 
 \int_0^t e^{- \beta s}dK_s^- \leq \int_0^t e^{- \beta s}dK'^-_s.$$
 \end{pro}
 
  \noindent{\textbf{Preuve. }} 1. 
   Suivant la construction de ces solutions donnée par \cite{martinLepl}, ces solutions 
  sont définies comme limite des  triplets  suivants: soient $(Y^n,Z^n,K^n)$ et $(Y'^{n},Z'^{n},K'^{n})$ solutions des EDSRR suivantes:
  
   {\footnotesize
      \begin{eqnarray} 
    Y^n_t & = & \int_t^T  e^{- \beta s}  f(s) ds+ \int_t^T  e^{- \beta s} dK_s^n -  n \int_t^T  e^{- \beta s}(Y^n_s-U_s)^+ds   -  \int_t^T  Z_s^n dW_s, ~ Y_T^{n} =0   \nonumber\\
     \label{ynt}\\
        Y^{'n}_t & = & \int_t^T  e^{- \beta s}  f'(s) ds+ \int_t^T  e^{- \beta s} dK_s^{'n}  -  n \int_t^T  e^{- \beta s}(Y^{'n}_s-U_s)^+ds -  \int_t^T  Z_s^{'n} dW_s, ~ Y_T^{'n} =0. \nonumber\\
   \end{eqnarray} 
   }
   On pose $f_n(t,y)=f(t) - ne^{- \beta t} (y-U_t)^+$ et $f'_n(t,y)=f'(t) - ne^{- \beta t} (y-U_t)^+$.
    Sous l'hypothèse  $(\mathcal{H}_1)$, nous obtenons    $f_n\leq f_n'$ et par suite, 
     pour tout $n$ (cf. proposition 2.3 \cite[p. 3]{hamadene}):
  $$ Y_t^n \leq Y_t'^n, ~    ~\forall ~ n,  t\geq 0, ~  ~ \mathbb{P}- \mbox{    p.s. et }
      \int_0^t  e^{- \beta s} dK_s^n\geq \int_0^t  e^{- \beta s} dK_s^{'n},~\forall t\leq n.$$   
      En faisant tendre $n$ vers l'infini,   puisque les suites de fonctions $f_n,f'_n$ sont décroissantes,     les limites décroissantes (respectivement croissantes) presque sûres 
  $ Y_t^n$ $Y_t'^n$ respectivement $\int_0^t  e^{- \beta s} dK_s^n$ et $ \int_0^t  e^{- \beta s} dK_s^{'n}$ 
    sont d'après respectivement les preuves du lemme 2  et du lemme 6 de \cite{martinLepl}
   $Y_t,~Y'_t,~\int_0^t e^{- \beta s}dK_s^{+}$ et  $\int_0^t e^{- \beta s}dK_s^{+'}$   
    qui vérifient bien à la limite
    $Y_t\leq Y'_t,$ $\int_0^t  e^{- \beta s} dK^+_s\geq \int_0^t  e^{- \beta s} dK'^+_s.$  \\
    \\
2.   Par définition des solutions des EDSRR en  horizon fini, il vient:  
$$\int_0^t  e^{- \beta s} dK_s^- =  Y_t - Y_0  +  \int_0^t  e^{- \beta s}  f(s) ds + \int_0^{t}  e^{- \beta s} dK_s^+   -  \int_0^{t}  Z_s dW_s.$$
De l'EDSRR (\ref{ynt}), on tire
$$ n \int_0^t  e^{- \beta s}(Y^n_s-U_s)^+ds   =  Y^n_t -Y^n_0 + \int_0^t  e^{- \beta s}  f(s) ds+ \int_0^t  e^{- \beta s} dK_s^n - \int_0^t  Z_s^n dW_s. $$
Par différence des deux précédentes égalités, on obtient: 
{\footnotesize $$n \int_0^t  e^{- \beta s}(Y^n_s-U_s)^+ds - \int_0^t  e^{- \beta s} dK_s^- = (Y^n_t-Y_t) - (Y^n_0 - Y_0) + \int_0^t  e^{- \beta s} (dK_s^n- dK_s^+)
- \int_0^t ( Z_s^n-Z_s ) dW_s.  $$}
Les suites  $(Y^n)$ et $(Z^n)$ sont deux suites de Cauchy  qui convergent
       respectivement dans $\mathbb{L}^2(\Omega)$ et $\mathbb{L}^2 ([0,\infty[ \times \Omega, dt \otimes d\mathbb{P})$ vers   $Y, Z$  (cf.  lemmes 5 et 6 de  \cite[p. 170-171]{martinLepl}). 
De plus, $\int_0^t  e^{- \beta s}  dK_s^n$ est une suite décroissante de limite presque sûre $\int_0^t  e^{- \beta s}  dK_s^+.$ \\
Par conséquent, $n \int_0^t  e^{- \beta s}(Y^n_s-U_s)^+ds $ converge p.s.  vers $\int_0^t  e^{- \beta s} dK_s^-$ lorsque $n$ tend vers l'infini.
\\
Ensuite, en utilisant $ Y_t^n \leq Y_t'^n$, nous obtenons:
$$ n \int_0^t  e^{- \beta s}(Y_s^n-U_s)^+ds   \leq n \int_0^t  e^{- \beta s}(Y_s'^n-U_s)^+ds. $$
Par passage à la limite, 
   nous obtenons la dernière comparaison recherchée: 
     $$ \int_0^t e^{- \beta s} dK_s^-\leq \int_0^t e^{- \beta s} dK'^-_s,~ ~ \mathbb{P}-\mbox{p.s.}$$ 
     \endproof

 \begin{lem}\label{lmmajU}
 Pour tout $n\geq0$ et $t\in [0,T],$  soit $(\overline{Y}^n,\overline{Z}^n,\overline{K}^n)$ solution de l'EDSRR à unique barrière   associée à 
$ (- \|f\|e^{-\beta t} - n e^{-\beta t} (y-U_t)^+,L)$ avec 
$ \overline{Y}^n_T =0,   \sup_t(L_t^+) \in \mathbb{L}^2$ et $ U_t = c_{1,0}e^{-\beta t}.$  Alors, 
\begin{equation}\label{majU}
\sup_{t \leq T}~  n(\overline{Y}^n_t - U_t)^+  \leq \beta  c_{1,0 } e^{-\beta t}, ~ ~ \mathbb{P}\mbox{-p.s.}
 \end{equation}
 
\end{lem}

 \noindent{\textbf{Preuve.}} C'est une adaptation de   la preuve  du lemme 41.4 de \cite[p. 165]{hamadene}.
  Pour tout $n,k\geq0,$ soit $(\overline{Y}^{n,k},\overline{Z}^{n,k})$ solution de l'EDS ordinaire suivante:
 {\footnotesize
 \begin{eqnarray*}
 \overline{Y}^{n,k}_t &  = &  -\int_t^T  e^{- \beta s} \|f\| ds  -  n \int_t^T  e^{- \beta s}(\overline{Y}^{n,k}_s-U_s)^+ds  
 + k \int_t^T  e^{- \beta s}(\overline{Y}^{n,k}_s-L_s)^-ds  -  \int_t^T  \overline{Z}^{n,k}_s dW_s  \nonumber\\
  \overline{Y}^{n,k}_T& =&0. 
 \end{eqnarray*}}
 Il est démontré dans \cite[p. 719]{Karoui1} que 
 $$\forall t \in [0,T], ~  \lim_{k \To \infty} k \int_t^T  e^{- \beta s}(\overline{Y}^{n,k}_s-L_s)^-ds = \int_t^T  e^{- \beta s} d\overline{K}^n_s, \quad \mathbb{P}\mbox{-p.s.}$$
et lorsque $k$ tend vers l'infini, nous avons 
 $$\lim_{k \To \infty} \overline{Y}^{n,k}_t  = \overline{Y}^{n}_t , \quad \mathbb{P}\mbox{-p.s.}$$ 
Ensuite, on définit $ X^{n,k}_t =\overline{Y}^{n,k}_t  -U_t, $ et par suite:
{\footnotesize  \begin{eqnarray*}
X^{n,k}_t &=& -U_t  -\int_t^T  e^{- \beta s} \|f\| ds  -  n \int_t^T  e^{- \beta s}(X^{n,k}_s)^+ds  
 + k \int_t^T  e^{- \beta s}(X^{n,k}_s-(L_s-U_s))^-ds  -  \int_t^T  \overline{Z}^{n,k}_s dW_s\nonumber\\
  X^{n,k}_T &=& - U_T. 
 \end{eqnarray*}}
  Appliquons la formule d'Itô entre $t$ et $T$ au produit $X^{n,k} \times \exp \big(- \int_0^. (\mu (u) + \nu (u))du \big),$ où 
  $\mu, \nu \in \mathcal{D}_n.$ Nous obtenons $\forall \nu, \forall \mu$:
{  \footnotesize
   \begin{eqnarray*} 
  X^{n,k}_t &=& 
   \EE \Big [-U_T \exp \big(- \int_t^T (\mu (u) + \nu (u))du \big)\nonumber\\
    &+& \int_t^T \exp \big(- \int_t^s (\mu (u) + \nu (u))du \big) (\dot{U}_s - e^{- \beta s} \|f\|  + \mu(s) (L_s-U_s)) ds\nonumber\\
    &+& \int_t^T \exp \big(- \int_t^s (\mu (u) + \nu (u))du \big)   \mu(s)(X^{n,k}_s-(L_s-U_s)) ds\nonumber\\
   & +& \int_t^T \exp \big(- \int_t^s (\mu (u) + \nu (u))du \big) \Big(ke^{- \beta s} (X^{n,k}_s-(L_s-U_s))^- + \nu(s)  X^{n,k}_s- n e^{- \beta s}( X^{n,k}_s)^+ \Big)ds \Big| \F_t\Big].
   \end{eqnarray*} 
 }
L'égalité $(6.13)$ de  Cvitanic et  Karatzas  \cite[p. 2042]{cvitanic } montre que 
 \begin{eqnarray*} 
   X^{n,k}_t &\leq & \ESS{\mu  \in \mathcal{D}_n}~ \INF{  \nu \in \mathcal{D}_n}
 ~    \EE \Big [-U_T \exp \big(- \int_t^T (\mu (u) + \nu (u))du \big)\nonumber\\
  &+ & \int_t^T \exp \big(- \int_t^s (\mu (u) + \nu (u))du \big) (\dot{U}_s - e^{- \beta s} \|f\|  + \mu(s) (L_s-U_s)) ds 
   \Big| \F_t\Big].
   \end{eqnarray*} 
Du fait que  $L_t \leq  U_t, U_T\geq 0$ et $\mu \geq 0,$ il vient: 
\begin{eqnarray*} 
  X^{n,k}_t  &\leq & \ESS{\mu  \in \mathcal{D}_n}~ \INF{  \nu \in \mathcal{D}_n}
 ~   \EE \Big [  \int_t^T \exp \big(- \int_t^s (\mu (u) + \nu (u))du \big) |\dot{U}_s |ds 
   \Big | \F_t\Big]\nonumber\\
  &\leq & \ESS{\mu  \in \mathcal{D}_n} 
 ~   \EE \Big [ 
      \int_t^T \exp \big(- \int_t^s (\mu (u) + n)du \big) |\dot{U}_s |ds 
   \Big | \F_t\Big]\nonumber\\
  &\leq & \EE \Big [  \int_t^T \exp \big(-     n(s-t) \big) |\dot{U}_s |ds 
   \Big | \F_t\Big], 
         \end{eqnarray*}  
les deux  dernières inégalités sont obtenues du fait que   $e^{-\nu} \geq e^{-n}$ et $e^{-\mu} \leq 1.$\\
  Puisque  $U_t = c_{1,0} e^{- \beta t},$ nous obtenons: 
$$X^{n,k}_t \vee 0 \leq    \frac{\beta   c_{1,0}}{ n } e^{- \beta t} ,$$ 
et  pour tout $t \in [0,T],$ il vient:   
$$
  n(\overline{Y}^n_t - U_t)^+  = \lim_{k \To \infty} n (X^{n,k}_t \vee 0 )\leq  \beta c_{1,0} e^{-\beta t}, ~ ~ \mathbb{P}\mbox{-p.s.}
 $$
 
\endproof

\begin{pro}\label{procompKPI}
Pour tout $n\geq0$ et $t\in [0,T],$ soient $(Y^n,Z^n,K^n)$  et $(\rho, \theta, \Pi)$ solutions de l'EDSR réfléchie à unique barrière  associées respectivement à $ (f e^{-\beta t} - n e^{-\beta t} (\overline{Y}^n_t-U_t)^+,L)$
  et $ (- \|f\|e^{-\beta t} - \beta c_{1,0}e^{-\beta t},L)$ avec $ Y^n_T =\rho_T =0, \sup_t(L_t^+) \in \mathbb{L}^2,  U_t = c_{1,0} e^{-\beta t}$
  et $\overline{Y}^n$ est défini dans le lemme \ref{lmmajU}. 
    Alors, 
$$ \forall t \in [0,T], ~ ~ Y^n_t  \geq \rho_t, ~  ~ \mbox{et}  ~    \int_0^t e^{-\beta s} dK^n_s  \leq  \int_0^t e^{-\beta s} d\Pi_s.$$
 De plus, pour $0 < \eps < \frac{1}{6 \exp(\beta^{-1})},$ il existe une constante $\tilde{C}_\eps$ telle que $\forall T:$
\begin{equation}\label{compKPI}
 \EE \Big (\int_0^T e^{-\beta s} d \Pi_s \Big  )^2 \leq \tilde{C}_\eps,
\end{equation}
où $$\tilde{C}_\eps = \Big (\frac{1}{3}-2\eps \exp\Big( \frac{1}{\beta }  \Big) \Big)^{-1}
  \left[  \Big(1+  2  \exp\Big( \frac{1}{\beta }\Big)\Big) \frac{1}{\beta} (\|f\|    +   \beta   c_{1,0}   )^2
  +\frac{1}{\eps}\EE[\sup_{s\geq 0}(L_s^+)^2]\right].$$

\end{pro}

  \noindent{\textbf{Preuve.}}  
  Pour tout $n\geq0$ et $t\in [0,T],$ soient $( Y^n, Z^n,K^n)$ et $(\rho, \theta, \Pi)$ solutions de l'EDSR associées respectivement à  à $ (f e^{-\beta t} - n e^{-\beta t} (\overline{Y}^n_t-U_t)^+,L)$
 et $ (- \|f\|e^{-\beta t} - \beta c_{1,0}e^{-\beta t},L)$ 
  avec $  Y^n_T =\rho_T =0,  \sup_t(L_t^+) \in \mathbb{L}^2$, $ U_t = c_{1,0} e^{-\beta t}:$
   \begin{eqnarray*}  
Y^n_t  & =&   \int_t^T  e^{- \beta s}  f(s) ds  -  n \int_t^T  e^{- \beta s}(\overline{Y}^n_s-U_s)^+ds  
 + \int_t^T e^{-\beta s} dK^n_s   -  \int_t^T  Z_s^n dW_s \nonumber\\
  \rho_t  & =&    - \int_t^T  e^{- \beta s}  \|f\|  ds  -    \int_t^T \beta c_{1,0}  e^{- \beta s} ds  
 +   \int_t^T e^{-\beta s} d\Pi_s -  \int_t^T  \theta_s dW_s. 
\end{eqnarray*}  
La   majoration (\ref{majU}) du  lemme \ref{lmmajU}  (avec $ \sup_t(L_t^+) \in \mathbb{L}^2$ et $ U_t = c_{1,0} e^{-\beta t}$) implique 
$$-  n(\overline{Y}^n_t - U_t)^+  \geq  -  \beta  c_{1,0 } e^{-\beta t}, ~ ~ \mathbb{P}\mbox{-p.s.}
 $$
  et naturellement $e^{- \beta s}f(s) \geq -e^{- \beta s}\|f\|. $
Par conséquent, 
$$  e^{- \beta s}  f(s) - n   e^{- \beta s}(\overline{Y}^n_s-U_s)^+  \geq 
-    e^{- \beta s}  \|f\|    -      \beta c_{1,0}  e^{- \beta s}, $$
et en utilisant
le corollaire de comparaison \ref{comp}, il vient: 

$$\forall t \in [0,T], ~ Y^n_t \geq  \rho_t, ~  ~ \mbox{et}  ~    \int_0^t e^{-\beta s} dK^n_s  \leq  \int_0^t e^{-\beta s} d\Pi_s ,~ ~ \mathbb{P}\mbox{-p.s.}$$
Cherchons maintenant une estimation de  $\EE\big (\int_0^T e^{-\beta s} d \Pi_s\big  )^2. $
Appliquons la formule d'Itô à $\rho^2_t$ entre $t$ et $T$: 
 $$\EE(\rho_t )^2 +   \EE \int_t^T  \theta_s^2 ds=     2\EE  \int_t^T  e^{- \beta s}  (-\|f\|    -     \beta c_{1,0} )\rho_s ds  + 2 \EE\int_t^T e^{-\beta s} 
\rho_s d\Pi_s.$$  
L'inégalité de Cauchy-Schwartz et  l'égalité $\int e^{-\beta s} (\rho_s - L_s) d\Pi_s =0$ impliquent  
$$
\EE(\rho_t )^2 +   \EE \int_t^T  \theta_s^2 ds \leq    \EE  \int_t^T  e^{- \beta s} \rho_s^2 ds +       \EE  \int_t^T  e^{- \beta s}  (\|f\|    +   \beta   c_{1,0}   )^2 ds  + 2 \EE\int_t^T e^{-\beta s} 
L_s d\Pi_s, 
$$
   de plus,   pour tout $\eps>0:$
\begin{eqnarray*}
2\int_t^{T} L_s e^{- \beta s} d\Pi_s^+ &\leq& 2\int_t^{T} L_s^+ e^{- \beta s} d\Pi_s^+\nonumber\\
&\leq&
   2\sup_{s\geq t}L_s^+\int_t^{T} e^{- \beta s} d\Pi_s^+\leq 
   \frac{1}{\eps}sup_{s\geq t}(L_s^+)^2+
  \eps\Big(\int_t^T e^{- \beta s}d\Pi_s^+\Big)^2
  \end{eqnarray*}
   (on utilise $2ab\leq \frac{1}{\eps}a^2+\eps b^2$).
    On applique le lemme de Gronwall     avec  
  $\psi (t) =  e^{- \beta t}$ et 
$$
\varphi(t) =  \frac{1}{\beta}    (\|f\|    +   \beta   c_{1,0}   )^2 e^{- \beta t} 
  +\frac{1}{\eps}\EE[sup_{s\geq t}(L_s^+)^2]+  \eps \EE\Big[\Big(\int_t^T e^{- \beta s}d\Pi_s^+\Big)^2\Big].
  $$
   Comme   $\varphi$ est décroissante   (cf. appendice) 
 il vient:
   \begin{equation}\label{estMM1}
   \EE[\rho_t^2]\leq \varphi(t) \exp \Big(\int_t^T  e^{- \beta u}  du\Big ) \leq  \exp\Big( \frac{1}{\beta }  \Big)\varphi(t).
    \end{equation}
  On a aussi par conséquent, puisque $\int_t^T|\theta_s|^2ds\leq \rho_t^2+\int_t^T|\theta_s|^2ds,$
  la même majoration
  \begin{equation}\label{estMM11}
  \EE(\int_t^T|\theta_s|^2ds)\leq  \exp\Big( \frac{1}{\beta }  \Big)\varphi(t).
   \end{equation}
  Ensuite, nous avons: 
$$ \int_0^T e^{-\beta s} d\Pi_s =   \rho_0   +  \int_0^T  e^{- \beta s}  \|f\|  ds  +    \int_0^T \beta c_{1,0}  e^{- \beta s} ds  
 +      \int_0^T  \theta_s dW_s. 
$$
 Des inégalités   (\ref{estMM1}) et (\ref{estMM11}), nous obtenons: 
  $$
 \frac{1}{3} \EE\big(\int_0^T e^{-\beta s} d\Pi_s \big)^2  
   \leq   \frac{1}{\beta} (\|f\|    +   \beta   c_{1,0}   )^2   + 2   \exp\Big( \frac{1}{\beta }  \Big)\varphi(0).  
 $$
  Si on retranche de $\varphi(0)$  le terme $\eps \EE\Big[\Big(\int_0^T e^{- \beta s}d\Pi_s^+\Big)^2\Big]$
    il vient
   \begin{equation}
   \label{majIntK1}
  \Big (\frac{1}{3}-2\eps \exp\Big( \frac{1}{\beta }  \Big) \Big)\EE\Big(\int_0^T e^{-\beta s} d\Pi_s \Big)^2  \leq  
    \Big(1+  2  \exp\Big( \frac{1}{\beta }\Big)\Big) \frac{1}{\beta} (\|f\|    +   \beta   c_{1,0}   )^2
  +\frac{1}{\eps}\EE[\sup_{s\geq 0}(L_s^+)^2],
   \end{equation}
 avec  $\eps$ est choisi assez petit vérifiant $\eps < \frac{1}{6 \exp(\beta^{-1})}.$
 Enfin,  $\forall T$:
$$
 \EE\Big(\int_0^T e^{-\beta s} d\Pi_s \Big)^2  \leq   \Big (\frac{1}{3}-2\eps \exp\Big( \frac{1}{\beta }  \Big) \Big)^{-1}
  \left[  \Big(1+  2  \exp\Big( \frac{1}{\beta }\Big)\Big) \frac{1}{\beta} (\|f\|    +   \beta   c_{1,0}   )^2
  +\frac{1}{\eps}\EE[\sup_{s\geq 0}(L_s^+)^2 ] \right].
$$

  \endproof

 \noindent{\textbf{Preuve du théorème \ref{bar2}.}}
  En utilisant le théorème 3.2 de \cite{hamadene}, on déduit l'existence  pour tout $n$
   d'un quadruplet de processus   qui vérifient $Y^n \in \mathcal{C}^2([0,n]),~  Z^n \in  \mathbb{H}^2 ([0,n] \times \Omega)$  et 
    $ ~ \forall t \leq n:$  
         
         $$ (\mathcal{S_{L,U}}) \left\lbrace
\begin{array}{lcl}
&&Y^n_t =0 ~ \forall t \geq n ,\\
 &&\forall t\leq n, L_t \leq Y^n_t \leq U_t, ~ 
 \\
 && Y^n_t   =   \int_t^n e^{- \beta s}  f(s) ds +  \int_t^n e^{- \beta s}  dK^{n+}_s  - \int_t^n e^{- \beta s}  dK^{n-}_s   -  \int_t^{n}  Z_s^n dW_s. \\
  \\
 &&  \mathbb{E} (\int_0^t   e^{- \beta s}  dK_s^{n+} )^2 < \infty ~ \mathbb{E} (\int_0^t e^{- \beta s}  dK_s^{n-} )^2 < \infty ~  \mbox{et}\\
 \\
  &&  \int_0^t  e^{- \beta s}  (Y_s^n - L_s)dK_s^{n+}  = \int_0^t  e^{- \beta s}  (U_s - Y_s^n  )dK_s^{n-} =  0, ~ ~t \leq n ~  \mathbb{P}\mbox{-p.s.}
\end{array}\right.
$$    
 
 La proposition \ref{comp2bar} implique que $(Y^n)_n$ est une suite croissante de processus,  $(\int e^{- \beta s} dK^{n+} )_n$ est une suite  décroissante  de processus et    $(\int e^{- \beta s} dK^{n-} )_n$  est une suite   croissante  de processus. On note leurs limites presque sûres respectivement $Y$, $\int e^{- \beta s} dK_s^+$   et $\int e^{- \beta s} dK_s^-. $
  Ensuite, $\forall k  \geq n$ nous obtenons:
  {\footnotesize
  \begin{eqnarray*}
 Y^{n+k}_t - Y^n_t   &=&    \int_n^{n+k} e^{- \beta s}  f(s)    ds 
 +   \int_t^n e^{- \beta s} \left[ dK^{(n+k)+}_s  -dK^{n+}_s  \right  ]   + \int_n^{n+k} e^{- \beta s}dK^{(n+k)+}_s \nonumber\\
  & - &  \int_t^n e^{- \beta s} \left[ dK^{(n+k)-}_s  -dK^{n-}_s  \right  ]  - \int_n^{n+k} e^{- \beta s}dK^{(n+k)-}_s \nonumber\\
 & - & \int_t^{n}\left[  Z_s^{n+k} - Z_s^{n} \right]  dW_s - \int_n^{n+k}  Z_s^{n+k} dW_s.
     \end{eqnarray*}
     }
  Appliquons la formule d'Itô à $( Y^{n+k}_t - Y^n_t  )^2$ entre $t$ et $n$, nous obtenons:
     \begin{eqnarray*}
   ( Y^{n+k}_t - Y^n_t  )^2  + \int_t^{n}\left(  Z_s^{n+k} - Z_s^{n} \right  )^2  ds 
    &=&   2\int_t^{n}e^{- \beta s}( Y^{n+k}_s - Y^n_s  ) \left(  dK_s^{(n+k)+} - dK_s^{n+} \right  )    
  \nonumber\\
   &-&   2\int_t^{n}e^{- \beta s}( Y^{n+k}_s - Y^n_s  ) \left(  dK_s^{(n+k)-} - dK_s^{n-} \right  )   
   \nonumber\\
 & - &  2\int_t^{n}( Y^{n+k}_s - Y^n_s  ) \left(  Z_s^{n+k} - Z_s^{n} \right  )  dW_s  +  (Y_n^{n+k})^2. 
     \end{eqnarray*}
  En utilisant la décroissance   de    la suite $(\int e^{- \beta s} dK^{n+
  } )_n$ et   la croissance de la suite  $(\int e^{- \beta s} dK^{n-} )_n$, puis 
    en passant à l'espérance, nous obtenons:
  $$ \mathbb{E}( Y^{n+k}_t - Y^n_t  )^2  + \mathbb{E} \int_t^{n}\left(  Z_s^{n+k} - Z_s^{n} \right  )^2  ds  \leq \EE(Y_n^{n+k})^2.$$
     Ensuite, en utilisant (\ref{estY_t3}) pour $L\leq 0 \leq U$ et  tout $n$ et tout $k$ 
   $$\exists D>0~~\EE(Y^{n+k}_s)^2\leq D e^{- \beta s}$$
   avec $D =  \frac{1}{\beta} \| f\|^2   \exp\Big( \frac{2 C+1}{\beta } \Big ),$
    et alors on a simplement
    $$ \mathbb{E}( Y^{n+k}_t - Y^n_t  )^2 
 + \mathbb{E} \int_t^{n}\left(  Z_s^{n+k} - Z_s^{n} \right  )^2  ds \leq \EE(Y^{n+k}_n)^2 \leq De^{- \beta n}\rightarrow 0.$$
     Par conséquent,    les suites  $(Y^n)$ et $(Z^n)$ sont deux suites de Cauchy qui convergent
       respectivement dans $\mathbb{L}^2(\Omega)$ et $\mathbb{L}^2 ([0,\infty[ \times \Omega, dt \otimes d\mathbb{P})$ vers deux processus $Y, Z:$  
    $$\lim_{n \To \infty}Y^n_t = Y_t ~ ~ \mbox{et} ~ ~ \lim_{n \To \infty}Z^n_t = Z_t.$$
    On a donc $Y_t\in \mathbb{L}^2$ pour tout $t$ et $Z\in \mathbb{L}^2 ([0,\infty[ \times \Omega, dt \otimes d\mathbb{P})$.\\
    \\
   $ii/$ Par limite presque sûre  
   de la relation  vérifiée pour tout $n$ par le quadruplet
   $(Y^n,Z^n,K^{n+},K^{n-})$, on obtient l'item $ii/.$
   \\
   De plus, puisque pour tout $n$ $L_t\leq Y^n_t\leq U_t$ à la limite on a également $L_t\leq Y_t\leq U_t$
   presque sûrement soit l'item $iii/$.\\
    \\
      $iv/$      La majoration (\ref{compKPI}) de la proposition \ref{procompKPI}   implique 
       $$\forall t, \forall t \leq n, ~ \EE \big( \int_0^t  e^{- \beta s}   dK^{n+}_s \big)^2 \leq\EE \big( \int_0^t  e^{- \beta s}   dK^{n,k}_s \big)^2 \leq \EE \big( \int_0^t  e^{- \beta s}   d\Pi_s \big)^2 \leq\tilde{C}_\eps$$
où $\int_0^t  e^{- \beta s}   dK^{n+}_s  = \displaystyle \lim_{k \To \infty }\int_0^t  e^{- \beta s}   dK^{n,k}_s$ et 
$$ \tilde{C}_\eps = \Big (\frac{1}{3}-2\eps \exp\Big( \frac{1}{\beta }  \Big) \Big)^{-1}
  \left[  \Big(1+  2  \exp\Big( \frac{1}{\beta }\Big)\Big) \frac{1}{\beta} (\|f\|    +   \beta   c_{1,0}   )^2
  +\frac{1}{\eps}\EE[\sup_{s\geq 0}(L_s^+)^2]\right].$$                  
     Quand $n$ tend vers l'infini, il vient: 
             $$\forall t, ~  \EE \big( \int_0^t  e^{- \beta s}   dK^+_s  \big)^2 \leq \tilde{C}_\eps.$$
             Puis, on fait tendre  $t$ vers l'infini et on obtient: 
              $$ \EE \big( \int_0^\infty  e^{- \beta s}   dK^+_s  \big)^2 \leq \tilde{C}_\eps,$$
et par suite $\int_0^\infty  e^{- \beta s}   dK^+_s \in \mathbb{L}^2.$   
 Ensuite,  on tire de (\ref{EDSR2}): 
    $$\int_0^\infty e^{- \beta s} dK_s^-  =  Y_t - Y_0 +    \int_0^\infty   e^{- \beta s}  f(s) ds + \int_0^\infty   e^{- \beta s} dK_s^+  -  \int_0^\infty   Z_s dW_s, $$ 
    qui   est la somme de  cinq termes de $\mathbb{L}^2$, et par conséquent $\int_0^\infty  e^{- \beta s}   dK^-_s  $ est aussi de carré intégrable.\\

   Pour achever la preuve de l'item $iv/$, pour tout $n$ on a l'identité presque sûre
 $$\int_0^t  e^{-\beta s}(Y_s^n-L_s)dK_s^{n+}=  \int_0^t  e^{-\beta s}(U_s- Y_s^n)dK_s^{n-}=    0,  ~  ~ t \leq n.$$ 
    La monotonie des suites de mesures $K^{n+}$ et  $K^{n-}$ montre que $\lim_n~p.s.\int_0^t e^{-\beta s}L_sdK_s^{n+}=
\int_0^t e^{-\beta s}L_sdK_s^+ $ et $\lim_n~p.s.  \int_0^t e^{-\beta s}U_sdK_s^{n-}=
\int_0^t e^{-\beta s}U_sdK_s^-. $
    Par ailleurs développons  (en rappelant que par définition $Y_s^n=0$ pour tout $s\geq n$) la différence suivante:
    $$\int_0^te^{-\beta s}Y_s^ndK_s^{n+}-\int_0^t e^{-\beta s}Y_sdK_s^{+}=$$
    $$\int_0^te^{-\beta s}Y_s^n(dK_s^{n+} -dK_s^{+})+\int_0^t e^{-\beta s}(Y^n_s-Y_s)dK_s^{+}.$$
    Le deuxième terme tend presque sûrement vers $0$ par convergence de Lebesgue monotone.
    Pour le premier, notons que $L_s\leq Y_s^n\leq Y_s$, il converge vers 0 par convergence de Lebesgue
    majorée. De façon similaire, 
    développons   la différence suivante:
    $$\int_0^te^{-\beta s}Y_s^ndK_s^{n-}-\int_0^t e^{-\beta s}Y_sdK_s^{-}=$$
    $$\int_0^te^{-\beta s}Y_s^n(dK_s^{n-} -dK_s^{-})+\int_0^t e^{-\beta s}(Y^n_s-Y_s)dK_s^{-}.$$
    Le deuxième terme tend presque sûrement vers $0$ par convergence de Lebesgue monotone. 
      Pour le premier, notons que $L_s\leq Y_s^n\leq Y_s \leq U_s$, il  
     converge vers 0 par convergence de Lebesgue majorée. Par suite, 
     $$\int_0^t  e^{-\beta s}(Y_s-L_s)dK_s^{+}=  \int_0^t  e^{-\beta s}(U_s- Y_s)dK_s^{-}=    0. $$
           Suite de l'item $i/$: De la relation (\ref{EDSR2})
    on tire:

 $$      Y_t = Y_0- \int_0^t e^{- \beta s}  f(s)  ds  + \int_0^t  e^{- \beta s} dK_s^+ -\int_0^t  e^{- \beta s} dK_s^- +  \int_0^t Z_s dW_s .$$
 On a donc puisqu'ici $f$ ne dépend pas de $y$ et dans ce cas est bornée: 

   $$\sup_{  t}|Y_t|\leq |Y_0| +\frac{ \|f\|}{\beta }  +\int_0^\infty e^{- \beta s} dK_s^+ + \int_0^\infty e^{- \beta s} dK_s^- + \sup_{  t}| \int_0^{\infty} Z_s dW_s|, $$
 c'est à dire une majoration par la somme de cinq termes dans $\LL^2: \EE[\sup_{  t}|Y_t|^2] <\infty$ et $Y \in \mathcal{C}^2.$

\endproof

   \section{Existence de $(Y^1,Y^2)$} 
    \label{existe}

     En utilisant  le    théorème    \ref{bar2} avec 
     $$f(s) = f(0,X_s^0) -f(1,X_s^1), ~ ~ L_t = - c_{0,1}  e^{- \beta t}  \leq 0 \leq 
        U_t= e^{- \beta t}  c_{1,0}, $$ 
         il existe un quadruplet de processus mesurables $(Y,Z,K^+,K^-)  $ tel que:
      {\footnotesize       
         $$
         (S)\left\lbrace
\begin{array}{lcl}
 &&Y \in \mathcal{C}^2,  Z\in \mathbb{H}^2 . \\
 \\
&&  Y_t =   \int_t^{+\infty}  e^{- \beta s} ( f(0,X_s^0) -f(1,X_s^1) )ds + \int_t^{+\infty}  e^{- \beta s} dK_s^+ - \int_t^{+\infty}  e^{- \beta s} dK_s^-  -  \int_t^{+\infty}  Z_s dW_s \\
  \\
 &&- c_{0,1}  e^{- \beta t}  \leq Y_t  \leq   e^{- \beta t}  c_{1,0}\\
 \\
 &&   \int_0^{+\infty}  e^{- \beta s} dK_s^{\pm} \in \mathbb{L}^2 ~ \mbox{et} ~ e^{- \beta t }(Y_t + e^{- \beta t}  c_{0,1} )dK_t^+  =   e^{- \beta t }(   c_{1,0}e^{- \beta t}   -Y_t )dK_t^-= 0.
\end{array}\right.$$    
}                      
 Nous pouvons donc prouver notre résultat principal: l'existence des processus $(Y^1,Y^2)$ introduits  dans la proposition \ref{defopt}. On se réfère au théorème 3.2  \cite[p. 186]{jeanblanc}.

     \begin{thm} Il existe un couple de processus continus $(Y^1_t,Y^2_t)_{t\geq 0}$  à valeurs dans $\R$ qui satisfont (\ref{Y1}) et (\ref{Y2}): 
     
     \begin{eqnarray*}
     Y_t^1 &=& \ESS{\theta \in \mathcal{T}_t}  ~  \mathbb{E} \left [ \int_t^{\theta}e^{- \beta s} f(0,X_s^0) ~ds   -  e^{- \beta \theta} c_{0,1} + Y_{\theta}^2 | \F_t\right],~Y^1_\infty=0   \\
       Y_t^2 &=& \ESS{\theta \in\mathcal{T}_t}  ~  \mathbb{E} \left [ \int_t^{\theta}e^{- \beta s} f(1,X_s^1) ~ds    - e^{- \beta \theta} c_{1,0}  + Y_{\theta}^1 | \F_t\right],~Y^2_\infty=0.
      \end{eqnarray*}

     \end{thm}

       \noindent{\textbf{Preuve. }}   On applique le théorème  \ref{bar2} avec $$f(s) = f(0,X_s^0) -f(1,X_s^1), ~ ~ L_t = - c_{0,1}  e^{- \beta t}  \leq 0 \leq         U_t= e^{- \beta t}  c_{1,0}, $$
      Puisque $\int_t^{+\infty}  e^{- \beta s} dK_s^{\pm}$ sont intégrables et $f$ bornée on
         peut  définir les processus 
             \begin{eqnarray*}
Y_t^1 & = &      \mathbb{E} \left[ \int_t^{ \infty}  e^{- \beta s}  f(0,X_s^0) ds   +  \int_t^{+\infty}  e^{- \beta s} dK_s^+ \big | \F_t  \right]\nonumber
\\
Y_t^2 &=& \mathbb{E} \left[ \int_t^{ \infty}  e^{- \beta s}  f(1,X_s^1) ds   +  \int_t^{+\infty}  e^{- \beta s} dK_s^- \big | \F_t  \right].
   \end{eqnarray*}
   Par suite, $  Y_t =   Y_t^1 -   Y_t^2, ~t \geq 0 $ et  de l'inégalité $- c_{0,1}  e^{- \beta t}  \leq Y_t  \leq   e^{- \beta t}  c_{1,0},$ on obtient $Y_{\infty} = 0.$   
   
 Montrons maintenant  que $\EE  [\sup_{t\geq 0}|Y_t^i|^2 ]< \infty, ~ i = 0,1.$ 
   En passant au carré et à l'espérance des précédentes égalités, en utilisant la bornitude de la fonction $f$,  et le fait que  $\int_0^t  e^{- \beta s} dK_s^+$ est positive, nous obtenons pour tout $t$:
 $$Y_t^1  \leq       \frac{1}{\beta}\|f\|  +  \EE[( \int_0^{+\infty}  e^{- \beta s} dK_s^+ )/\F_t].$$
 L'inégalité de Burkholder-Davis-Gundy  \cite[p. 166]{karatzas} 
 appliquée à la martingale de carré intégrable $\EE[( \int_0^{+\infty}  e^{- \beta s} dK_s^+ )/\F_t]$ permet de conclure
 que $\EE  [\sup_{t\geq 0}|Y_t^1|^2 ]< \infty.$  
 \\
 Nous procédons  de la même façon pour démontrer  $\EE  [\sup_{t\geq 0}|Y_t^2|^2 ]< \infty.$ 
   \\
   
   Ensuite, d'après le théorème de représentation des martingales, il existe un processus prévisible adapté $Z^1$ 
    tel  que pour tout $t$:
    {\footnotesize 
   $$  \int_0^t  Z_s^1 dW_s= \int_0^t  e^{- \beta s}  f(0,X_s^0) ds +\int_0^t  e^{- \beta s} dK_s^+ + Y^1_t- \EE\Big[\int_0^\infty e^{- \beta s}  f(0,X_s^0) ds +\int_0^\infty e^{- \beta s} dK_s^+\Big] .$$}
En utilisant la troisième inégalité   du système   (S), nous avons $Y_t \geq -  c_{0,1} e^{- \beta t} $, et en remplaçant $Y_t$ par $Y_t^1 -Y_t^2$, nous obtenons $Y^1_t \geq -  c_{0,1} e^{- \beta t}+Y^2_t.$ \\
\\
De même,  la quatrième égalité   du système  (S), c'est à dire 
$\int_0^{.}  e^{- \beta t }(Y_t +  c_{0,1} e^{- \beta t})dK_t^+ =0,$ en  remplaçant $Y_t$ par $Y_t^1 -Y_t^2$ montre
 $$\forall T, ~ \int_0^{T}  (Y^1_t  - Y^2_t  +   c_{0,1}e^{- \beta t} )e^{- \beta t } dK^+_t =0.$$

Par suite, le triplet $(Y^1,Z^1,K^+)$ satisfait:
  $$
\left\{
\begin{array}{ccc}
  - dY^1_t = e^{- \beta t} f(0,X_t^0) dt +   e^{- \beta t} dK_t^+  -     Z_t^1 dW_t \\
Y^1_t \geq -  c_{0,1}e^{- \beta t } + Y^2_t\mbox{ et }e^{- \beta t }(Y^1_t- Y^2_t  +  e^{- \beta t } c_{0,1} ) dK^+_t =0.
\end{array}
\right.
$$
 Ensuite, on utilise le théorème \ref{bar1} item $v/$ avec $L_t  =- c_{0,1}e^{- \beta t} + Y^2_t.$
  Du fait que $\EE  [\sup_{t\geq 0}|Y_t^2|^2 ]< \infty,$ l'hypothèse $\EE  [\sup_{t\geq 0}(L_t^+)^2 ]< \infty$ est vérifiée et on a :
  $$ Y_t^1 = \ESS{\theta \in \TT} ~ \mathbb{E} \left[ \int_t^{ \theta}  e^{- \beta s}  f(0,X_s^0) ds   - c_{0,1} e^{- \beta \theta}+ Y^2_\theta \big | \F_t  \right], ~  t \geq 0.
$$
  D'une manière similaire, en utilisant la troisième inégalité   du système   (S), nous avons $Y_t \leq   c_{1,0}e^{- \beta t}$, et en remplaçant $Y_t$ par $Y_t^1 -Y_t^2$, nous obtenons $Y^2_t \geq -  c_{1,0}e^{- \beta t} + Y^1_t.$ 
Par suite,  le triplet $(Y^2,Z^1,K^-)$ satisfait:
  $$
\left\{
\begin{array}{ccc}
  - dY^2_t = e^{- \beta t} f(1,X_t^1) dt +   e^{- \beta t} dK_t^-  -     Z_t^1 dW_t \\
Y^2_t \geq -   c_{1,0} e^{- \beta t}+ Y^1_t  ~  \mbox{et} ~  e^{- \beta t}(  c_{1,0}e^{- \beta t} - Y^1_t  + Y^2_t ) dK^-_t =0.
\end{array}
\right.
$$
  Ensuite, on utilise le théorème \ref{bar1} item $v/$ avec $L_t  =- c_{1,0}e^{- \beta t} + Y^1_t.$
   Du fait que   $\EE  [\sup_{t\geq 0}|Y_t^1|^2 ]< \infty,$ l'hypothèse $\EE  [\sup_{t\geq 0}(L_t^+)^2 ]< \infty$ est vérifiée et on a :
  $$Y_t^2 = \ESS{\theta \in \TT} ~ \mathbb{E} \left[ \int_t^{ \theta}  e^{- \beta s}  f(1,X_s^1) ds   
  - c_{1,0}e^{- \beta \theta} + Y^1_\theta \big | \F_t  \right], ~  t \geq 0.
$$
D'où l'existence du couple $(Y^1,Y^2)$. 
 
     \endproof
    
        \section{Appendice}
    \textit{Lemme de Gronwall: } Si  $y$ et  $\psi$ sont deux fonctions continues   
    qui vérifient:
$$y(t)   \leq \varphi (t) + \int_t^T \psi (s)y(s) ds,$$
alors 
$$y(t)  \leq \varphi (t) + \int_t^T \varphi(s)  \psi (s) \exp\Big(\int_t^s  \psi (u) du\Big) ds.$$
Si de plus la fonction $\varphi$ est décroissante, alors: 
$$y(t)  \leq \varphi (t)  \exp\Big(\int_t^T  \psi (u) du\Big).$$


\end{document}